\documentclass{article}

\usepackage{microtype}
\usepackage{graphicx}
\usepackage{subcaption}
\usepackage{booktabs}
\usepackage{hyperref}

\usepackage[accepted]{icml2024}

\usepackage{amsmath}
\usepackage{amssymb}
\usepackage{mathtools}
\usepackage{amsthm}
\usepackage[capitalize,noabbrev]{cleveref}

\theoremstyle{plain}
\newtheorem{theorem}{Theorem}[section]

\newtheorem{lemma}[theorem]{Lemma}
\newtheorem{corollary}[theorem]{Corollary}
\theoremstyle{definition}

\theoremstyle{remark}

\usepackage[textsize=tiny]{todonotes}

\usepackage{times}
\usepackage{soul}
\usepackage{url}
\usepackage[utf8]{inputenc}
\usepackage{graphicx}
\usepackage{amsmath}
\usepackage{amsthm}
\usepackage{booktabs}
\usepackage{algorithm}
\usepackage{algorithmic}
\usepackage{amsfonts}       
\usepackage{xcolor}         
\usepackage{multirow}
\usepackage{bm}
\usepackage{pgfplotstable}
\usepgfplotslibrary{groupplots}
\usetikzlibrary{matrix}
\pgfplotsset{compat=newest}
\definecolor{marine}{RGB}{0,32,96}
\definecolor{navy}{RGB}{0,0,128}
\definecolor{maroon}{RGB}{128,0,0}
\definecolor{olivegreen}{RGB}{85,107,47}
\definecolor{gray}{RGB}{102,102,102}
\definecolor{green}{RGB}{131,198,210}
\definecolor{skyblue}{rgb}{0.3010, 0.7450, 0.9330}
\definecolor{purple}{rgb}{0.4940, 0.1840, 0.5560}
\definecolor{orange}{rgb}{0.9290, 0.6940, 0.1250}
\definecolor{brown}{RGB}{161,121,124}
\definecolor{deepblue}{rgb}{0.0, 0.0, 1.0}

\icmltitlerunning{Adaptive Stabilization based on Machine Learning for Column Generation}

\begin{document}

\twocolumn[
\icmltitle{Adaptive Stabilization Based on Machine Learning for Column Generation}




\begin{icmlauthorlist}
\icmlauthor{Yunzhuang Shen}{uts}
\icmlauthor{Yuan Sun}{ltu}
\icmlauthor{Xiaodong Li}{rmit}
\icmlauthor{Zhiguang Cao}{smu}
\icmlauthor{Andrew Eberhard}{rmit}
\icmlauthor{Guangquan Zhang}{uts}
\end{icmlauthorlist}

\icmlaffiliation{uts}{Australian Artificial Intelligence Institute, University of Technology Sydney, Australia}
\icmlaffiliation{ltu}{La Trobe Business School, La Trobe University, Australia}
\icmlaffiliation{rmit}{School of Computing Technologies, Royal Melbourne Institute of Technology, Australia}
\icmlaffiliation{smu}{School of Computing and Information Systems, Singapore Management University, Singapore}

\icmlcorrespondingauthor{Yunzhuang Shen}{shenyunzhuang@outlook.com}

\icmlkeywords{Machine Learning, ICML}

\vskip 0.3in
]



\printAffiliationsAndNotice{}  

\begin{abstract}

Column generation (CG) is a well-established method for solving large-scale linear programs. It involves iteratively optimizing a subproblem containing a subset of columns and using its dual solution to generate new columns with negative reduced costs. This process continues until the dual values converge to the optimal dual solution to the original problem. A natural phenomenon in CG is the heavy oscillation of the dual values during iterations, which can lead to a substantial slowdown in the convergence rate. \emph{Stabilization} techniques are devised to accelerate the convergence of dual values by using information beyond the state of the current subproblem. However, there remains a significant gap in obtaining more accurate dual values at an earlier stage. To further narrow this gap, this paper introduces a novel approach consisting of 1) a \emph{machine learning} approach for accurate prediction of optimal dual solutions and 2) an \emph{adaptive stabilization} technique that effectively capitalizes on accurate predictions. On the graph coloring problem, we show that our method achieves a significantly improved convergence rate compared to traditional methods. 
\end{abstract}

\section{Introduction}

Column generation (CG) is an effective method for solving linear programs (LP) with a large number of variables (or columns)~\citep{lubbecke2005selected}. It has many applications in solving combinatorial optimization problems with a decomposable structure~\citep{vanderbeck2000dantzig}, such as the vehicle routing problem~\citep{agarwal1989set}, the cutting stock problem~\citep{gilmore1961linear}, and the graph coloring problem~\citep{mehrotra1996column}.


CG solves a large-scale LP in iterative steps. In an iteration, a set of dual values is obtained by solving the LP that contains a small subset of columns and is then used to generate new columns with negative reduced costs. As this process repeats, the dual values converge to optimal dual values, i.e., an optimal dual solution to the original LP. This point of convergence is identified when no column with a negative reduced cost can be further generated.

\begin{figure}[tb!]
    \centering
    \includegraphics[width=0.99\columnwidth]{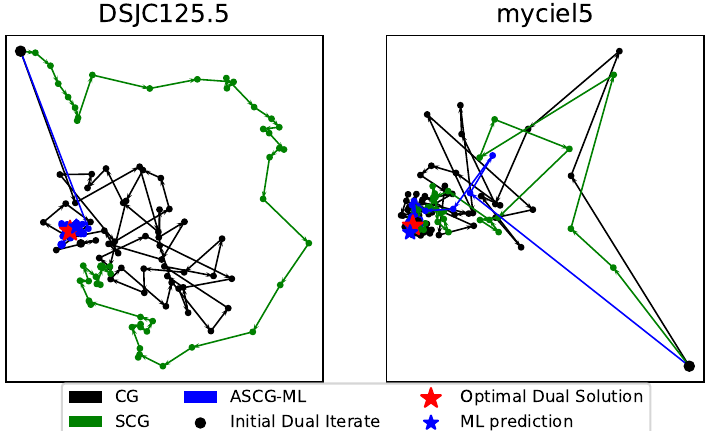}
        \caption{The trajectory of dual iterates for our adaptive stabilized CG based on ML (ASCG-ML), compared to CG and stabilized CG (SCG) for the first $50$ iterations. Each point (or a dual iterate) represents a set of dual values at an iteration of CG. The axes are the first two principal components of the high-dimensional dual space. The test graphs DSJC125.5 and myciel5 are from standard benchmarks for the graph coloring problem.}
        \label{fig:cg}
\end{figure}

As CG updates dual values by iteratively re-optimizing an evolving subproblem, this method may lead to significant oscillations in the dual iterates within the high-dimensional dual space. In this context, a dual iterate refers to the set of dual values obtained in a specific iteration of CG. This phenomenon is depicted in Figure~\ref{fig:cg}, which shows the trajectory of dual iterates for CG marked in green. Notably, the dual iterates tend to stay far from the optimal dual solution until a later stage. These issues can lead to the generation of redundant columns, causing a significant slowdown in the convergence rate.

Various techniques, termed \emph{stabilization}, have been devised to overcome this challenge. These techniques succeed in deriving dual values that are more closely aligned with the optimal dual solution by using information beyond the current state of the subproblem. Typically, this additional information includes the problem data~\citep{agarwal1989set, briant2008comparison, kraul2023machine} and/or historical dual iterates collected during the solution process~\citep{du1999stabilized, pessoa2018automation}. However, there is still significant room to obtain more accurate dual values at an earlier stage to further accelerate CG.

In this paper, we propose a new stabilization method to further narrow this gap. Our method consists of two parts: 1) a \emph{machine learning} (ML) approach to effectively predict the optimal dual solution, and 2) an \emph{adaptive stabilization} method that effectively uses ML prediction to guide the convergence of dual values. Specifically, we penalize dual variables taking values that deviate from ML prediction during subproblem optimization, resulting in a \emph{basin-in-attraction} effect pulling the dual values toward ML prediction. The strength of this attraction is adapted in line with the progression of CG, motivated by the fact that the accuracy of dual values can continuously improve and surpass the ML prediction. We name our method Adaptive Stabilized CG based on ML, ASCG-ML, and show that this method is grounded in Lagrangian duality theory and inherits the convergence guarantee from CG.

We empirically demonstrate the efficacy of ASCG-ML in the graph coloring problem (GCP)~\citep{jensen2011graph}. As illustrated in Figure~\ref{fig:cg}, ASCG-ML allows the dual values to be significantly drawn towards the prediction of the optimal dual solution, compared to standard CG and a stabilization approach~\citep{du1999stabilized}. This attraction is particularly strong when the dual values greatly differ from the optimal dual solution, as evidenced by a large `step size'. As CG progresses, our method consistently yields dual iterates that are in proximity to the optimal dual solutions, in contrast to the considerable fluctuations observed in the dual iterates produced by the compared methods. 

Our main contributions can be summarized as follows.

\begin{itemize}
    \item We formulate the task of identifying good dual values as a regression problem to predict optimal dual solutions. We develop two ML approaches to capture variable interrelationships for accurate predictions, including a feedforward neural network with statistical features and a graph neural network. 

    \item We propose ASCG-ML, which effectively uses an accurate ML prediction to guide the convergence of dual values. We theoretically show that ASCG-ML adaptively modulates the intensity of this guidance according to the progression of CG, and empirically show that ASCG-ML significantly improves the convergence rate compared to traditional approaches on the GCP. 
    
    \item We conduct an extensive comparative study of ASCG with two distinct ML approaches on both in-distribution and out-of-distribution data. Our findings present strong evidence that ASCG-ML can generalize across diverse, challenging families of benchmark instances beyond those it was trained on and pave the way for future advancements in the robust and generic integration of CG and ML techniques.
\end{itemize}

\section{Background \& Related Work}

\subsection{Column Generation for Graph Coloring}

CG, or Dantzig-Wolfe decomposition, is an effective method that can be used to obtain tight LP relaxation bounds for a wide range of combinatorial optimization problems~\citep{vanderbeck2000dantzig}. This is often crucial for exact solvers such as branch-and-bound~\citep{dakin1965tree}, the efficiency of which hinges on the strength of the bounding function. 

We will illustrate CG in the GCP, a classic NP-hard problem that finds many real-world applications~\citep{formanowicz2012survey} such as timetable scheduling. CG overcomes the highly symmetric characteristic of the GCP and manages to obtain tight LP relaxation bounds compared to alternative formulations of the problem. 

The Dantzig-Wolfe decomposition of the GCP can be expressed as using a minimum number of maximal independent sets~(MISs) to cover all vertices in a graph~\citep{mehrotra1996column}. Let $\mathcal{G}(\mathcal{V}, \mathcal{E})$ be a graph with the set of vertices $\mathcal{V}$ and the set of edges $\mathcal{E}$\footnote{Notations used in this paper are summarized in Appendix~\ref{appendix:notation}.}. Let $s \in \mathcal{S}$ be the set representation of a MIS, where $\mathcal{S}$ is the set of MISs in the graph. We describe the CG process in the dual LP relaxation to facilitate the later illustration of stabilization techniques, and note that the primal LP form is outlined in Appendix~\ref{appendix:gcp}. The dual LP relaxation of GCP can be defined as follows,
\begin{align}
    z^* = \max_{\bm{\pi}} \;& \sum_{i\in \mathcal{V}} \pi_i, & (\text{Dual Problem}) \label{eq:dualobj}\\
        s.t. \;& \sum_{i \in s} \pi_i \leq 1, &  s \in \mathcal{S},  \label{eq:dp-2}\\
               & 0 \le \pi_i \le 1, & i\in \mathcal{V}. \label{eq:dp-3}
\end{align}
Here, $\pi_i$ is a dual variable associated with the vertex $i$. A constraint in Eq.~\eqref{eq:dp-2} specifies that, for all vertices covered by this MIS (i.e., the expression $i \in s$), the sum of their corresponding dual values should not exceed $1$. Also note that each constraint in Eq.~\eqref{eq:dp-2} corresponds to a column in its primal form, and they represent the same MIS object. Since there may be exponentially many MISs in the graph, this LP typically cannot be solved directly using standard techniques such as simplex methods~\citep{dantzig2016linear}. 

CG tackles this large-scale dual problem in iterative steps. Specifically, a set of dual values $\bm{\pi}$ is obtained by solving the dual problem restricted to a partial set of MISs, referred to as the restricted dual problem (RDP). The dual values can be infeasible to the original dual problem, and are used to search for new MISs with minimum reduced costs. This process is called pricing and can be defined as the following combinatorial optimization problem,
\begin{align}
    c^* = \min_{\bm{v}} \;& 1 - \sum_{i \in \mathcal{V}} \pi_i\cdot v_i, & \text{(Pricing)}\\
        s.t. \;& v_i + v_j \leq 1, & (i,j)\in \mathcal{E} \\
              & v_i \in \{0,1\}, & i \in \mathcal{V},
\end{align}
where the binary decision variable $v_i$ denotes whether a vertex is in the solution, and solution space is the set of MISs. Note that this pricing problem involves solving an NP-hard maximum weight independent set problem. 

If the minimum reduced cost $c^*$ is negative, RDP is updated with the optimal MIS $s^*$ and re-optimized to produce new dual values with improved feasibility. Otherwise, the dual values are feasible for the dual problem, and they form an optimal dual solution. 


We have described CG in the context of GCP. However, this methodology is versatile and has proven effective in various problems~\citep{barnhart1998branch}. A generic description of CG involves iteratively improving the feasibility of dual values by generating columns with negative reduced costs. In different CG applications, their dual problems often bear resemblance to the structure presented in Eqs. \eqref{eq:dualobj}-\eqref{eq:dp-3}. However, a column can represent a different type of combinatorial object according to the structure of the pricing problem, such as a route in vehicle routing problems or a cutting pattern in the cutting stock problem.

\subsection{Stabilized Column Generation}
\label{subsec:scg}

A longstanding challenge in CG is the heavy oscillation of dual values. Specifically, the dual values are updated by the solution to the RDP. The evolution of the RDP causes changes in its associated dual polyhedron. Consequently, the solution to the RDP, which represents an extreme point of the dual polyhedron, can undergo significant shifts. This can lead to substantial oscillations in the dual iterates (Figure~\ref{fig:cg}) and a slowdown in the convergence rate.

Stabilization techniques aim to address this challenge and quickly produce dual values relevant to the optimal dual solution~\citep{lubbecke2005selected}. Several classes of stabilization methods have been proposed, including the \emph{limitation} approach, which imposes additional constraints on the RDP to restrict the dual space without cutting off all optimal dual solutions~\citep{marsten1975boxstep, ben2006dual, gschwind2017stabilized}; the \emph{centralization} approach, which generates columns using a central point of the dual polyhedron rather than a RDP solution~\citep{goffin1992decomposition,lee2011chebyshev}; the \emph{smoothing} approach, which blends the RDP solution with a reference point of good dual values~\citep{pessoa2018automation, wentges1997weighted}; and interior point methods with an early stopping criterion to solve the RDP~\citep{gondzio1996column, gondzio2013new}. These stabilization methods are not in direct competition with our proposed ML-based approach. For example, a ML prediction can be used as a reference point for the smooth approach, and our proposed method can be used in a complementary way to the limitation approach. 

Our method falls into the category of the \emph{penalization} approach. This approach was first proposed by \citet{du1999stabilized}, which penalizes dual variables taking values that deviate from the previous dual iterate. \citet{amor2009choice} studied the effect of different penalty functions. \citet{kraul2023machine} explored the prediction of optimal dual solutions through ML for the cutting stock problem and demonstrated that a feedforward neural network with basic problem attributes is sufficient to obtain high-quality prediction. In contrast, we study the dual solution prediction problem for the GCP and will show that modeling the interrelationships between variables is crucial for prediction accuracy. Moreover, we propose a novel adaptive stabilization method to optimize the utilization of ML prediction.




\subsection{Machine Learning For Optimization}

ML has been considered a promising technique in the broader context of optimization~\citep{ bengio2020machine}. Existing studies have shown that ML can enhance both exact methods~\citep{khalil2016learning, gasse2019exact, Morabit2020mlcs, babaki2021neural,deza2023machine} and heuristic methods~\citep{kool2018attention, sun2019using, xin2021neurolkh, liu2022learning, kotary2022fast, zong2022mapdp,huang2023searching,zhou2023towards}.

Among the exact approaches related to CG, several studies investigate the use of ML to select promising columns rather than the minimum-reduced-cost column selection rule~\citep{Morabit2020mlcs, babaki2021neural, chi2022deep}. Some others devise ML-enhanced heuristic pricing methods to efficiently solve pricing problems and generate high-quality columns~\citep{shen2022enhancing,quesnel2022deep}. Our work is complementary to the aforementioned studies, which require setting up the pricing problem with accurate dual values in the first place. 

Among ML-based heuristics, the most relevant to our study is a class of methods that develops an ML model to predict optimal solution values for integer variables, and uses ML prediction to improve heuristic search~\citep{li2018combinatorial, ding2020accelerating, shen2021learning,sun2022boosting}. Unlike this line of work, we study ML to make accurate predictions of the optimal dual solution in the context of CG. More specifically, the prediction targets are continuous variables, and our aim is to use ML prediction to accelerate the CG method to derive tight LP bounds.

\section{Approach}
\label{sec:ASCG}

In this section, we present ASCG-ML that can effectively derive dual values more relevant to the optimal dual solution compared to the RDP solution. Specifically, we consider a generalized formulation of the dual problem defined as
\begin{align}
    z^*_{\epsilon} = \max_{\bm{\pi}, \bm{w}^+, \bm{w}^-} & \sum_{s \in \mathcal{S}} (\pi_i - \epsilon  w_i^{-} - \epsilon w_i^{+}), \\
        s.t. & \sum_{i \in \mathcal{V}, i\in s} \pi_i \leq 1, &  s \in \mathcal{S},  \label{eq:dualcut}\\
             & -\pi_{i} - w_i^{-} \leq - \hat{y}_{i}, &  i \in \mathcal{V},  \label{eq:mdp-3}\\
             & \pi_{i} - w_i^{+} \leq \hat{y}_{i}, & i \in \mathcal{V},  \label{eq:mdp-4}\\
               & 0 \le \pi_i, w_{i}^{-}, w_{i}^{+} \le 1, & i\in \mathcal{V}. 
\end{align}
Compared with the standard dual problem in Eqs.~\eqref{eq:dualobj}-\eqref{eq:dp-3}, this generalized dual problem introduces new variables $w_i$ to measure the distance between the value taken by a dual variable $\pi_i$ and its predicted value by ML $\hat{y}_i$, as defined in the new constraints in Eqs.~\eqref{eq:mdp-3} and~\eqref{eq:mdp-4}. The parameter $\epsilon$ denotes the objective coefficient of the variable $w_i$, and the standard dual problem can be obtained by setting $\epsilon=0$ as a special case. With $\epsilon > 0$, it incurs a penalty with dual variables taking values deviating from the ML prediction, creating a {\em basin-in-attraction} effect from the ML prediction to the dual values. The magnitude of the penalty parameter $\epsilon$ controls the strength of this effect. 

We denote the generalized dual problem restricted to a subset of constraints by G-RDP. In contrast to the RDP, pricing using the solution to the G-RDP does not guarantee a nonpositive minimum reduced cost. This result is shown in Lemma~\ref{lemma:1} and will influence our design of the adaptive stabilization method in Section~\ref{subsec:apm}. All proofs can be found in Appendix~\ref{appendix:proof}.

\begin{lemma}
    \label{lemma:1}
    Let $\bm{\pi}^{\epsilon}$ denote the dual solution to the current G-RDP. The minimum reduced cost $c^*_{\epsilon}$ with respect to the dual values $\bm{\pi}^{\epsilon}$ can be positive. 
\end{lemma}

\begin{algorithm}[t!]\small
 \caption{ASCG-ML}
\label{alg:ada}
 \begin{algorithmic}[1]
    \REQUIRE An initial set of MISs $\mathcal{S}$, A ML prediction $\hat{\bm{y}}$;
    \STATE Initialize $\epsilon$;
    \REPEAT
        \STATE $\bm{\pi}^{\epsilon} \leftarrow$ solve G-RDP($\mathcal{S}, \hat{\bm{y}}, \epsilon$);
        \STATE $c^*, s^* \leftarrow$ pricing w.r.t. $\bm{\pi}^{\epsilon}$;
        \STATE $\mathcal{S} \leftarrow {s^*} \cup \mathcal{S}$;

        \STATE Update $\epsilon$ value;
    \UNTIL{$c^* \ge 0$ and $\epsilon = 0$} 
\end{algorithmic}
\end{algorithm}

Building on G-RDP, we outline the ASCG-ML method in Algorithm~\ref{alg:ada}. Two key components of ASCG-ML are the design of ML methods for accurate predictions and the adaptation of the penalty value during the iterative process, which are introduced in the following. 

\subsection{Dual Solution Prediction}\label{sec:DualPred}

We formulate the task of determining good dual values as a regression problem, where a ML model can be efficiently trained to minimize the mean squared error between the optimal solution values of the dual variables and the predicted values. We form the training data $\mathbb{D}^{(j)} = \{(\bm{x}^{i}, y^{i})| i \in \mathcal{V}\}$ from every optimally solved problem instance $j$. Each training example $(\bm{x}^{i}, y^{i})$ is associated with a dual variable $\pi_{i}$, where $\bm{x}^{i}$ denotes the feature vector of the variable and $y^{i}$ is its optimal solution value. Let $\hat{y}_{\theta} = f_{\theta}(\bm{x})$ denote the prediction of an ML model $f$ parameterized by $\theta$. Our loss function is defined as follows: 
\begin{equation}
\label{eq:reg}
    \mathcal{L}_{\theta} = \frac{1}{\sum_{j=1}^{J}|\mathbb{D}^{(j)}|} \sum_{j=1}^{J}\sum_{i=1}^{|\mathbb{D}^{(j)}|} (y^{j,i}-\hat{y}^{j,i}_{\theta})^{2}.
\end{equation}
We note that multiple optimal dual solutions may exist for an LP due to primal degeneracy, and when this occurs, we set the target value $y^{i}$ to an average of different optimal solution values of the corresponding variable. The benefit of this labeling approach will be examined empirically. 

We design two ML approaches for this regression problem, based on the understanding that modeling interrelationships between variables is critical for accurate predictions. For example, if a dual variable has a large value in the optimal dual solution, the variables in the same column with this variable must have relatively small optimal dual values. This is because the optimal dual values should satisfy all constraints in Eq.~\eqref{eq:dp-2}, i.e. the sum of the dual values in any column should not exceed $1$. It can also be noted that the importance of modeling variable relationships is recognized in a wider scope of ML applications for combinatorial optimization, such as predicting effective branching variables~\citep{gasse2019exact} in mixed integer programming solvers and predicting optimal integer solutions~\citep{ding2020accelerating}.

Our first approach involves equipping a feedforward neural network (FFNN) with statistical measures, each capturing a type of relation for a dual variable with others. These measures are computed over a set of randomly sampled columns and can be described using graph terminologies as follows: 1) the frequency of a vertex appearing in the samples, 2) statistics related to the cardinality of samples containing a vertex, 3) statistics of the average degree in samples containing a vertex, and 4) the degree of a vertex and graph density. The mathematical definition of these measures can be found in Appendix~\ref{appendix:mlmodel}.

Given this feature representation of a variable, a FFNN can be trained to predict its optimal solution value. A layer $l$ in our FFNN is defined as 
\begin{equation}
  \bm{h}^{(l+1)} = ReLU(W^{(l)}\bm{h}^{(l)} + \bm{b}^{(l)}),  
\end{equation}
where $W$ and $\bm{b}$ denote the learnable weight matrix and bias vector, and $\bm{h}$ is the hidden feature vector. ReLU is an activation function defined as $ReLU(\cdot) = max(0, \cdot)$.

Our second approach employs a graph convolutional network (GCN)~\citep{kipf2016semi} that automatically learns the interrelationships between variables and predicts optimal dual solutions. Specifically, in GCP, a problem instance corresponds to a graph and each dual variable corresponds to a vertex. Therefore, the graph naturally represents the relationships among the variables. Given this graph, the GCN predicts the optimal dual solution for the associated problem instance. A layer in this GCN $l$ is defined as 
\begin{equation}
  H^{(l+1)} = ReLU(\Tilde{D}^{-\frac{1}{2}}\Tilde{A}\Tilde{D}^{-\frac{1}{2}}H^{(l)}W^{(l)} + H^{(l)}),  
\end{equation}
where $\Tilde{A} = A + I$ is the adjacency matrix of the undirected graph with added self-connections, $\Tilde{D}_{ii} = \sum_i \Tilde{A}_{ij}$, and $H$ is a matrix composing hidden feature vectors of dual variables. The second term inside the activation function is commonly called the residual connection, and its purpose is to enable the training of deep neural networks with many layers. In the last layer, a linear mapping is used for both neural networks to output a scalar value for each variable as the prediction.

We provide a toy example to illustrate the intuition of GCN automatically learning useful information in our context. Consider that the feature vector of a variable (or vertex) has two dimensions. The first layer of GCN updates the first entry in the feature vector with the number of $1$-hop neighbors, and the second layer updates the second entry with the number of $2$-hop neighbors. If the number of $1$-hop neighbors for a variable is large, the variable tends to have a large dual solution value. Consider the extreme case where the variable is connected to all other variables; then this variable itself forms a maximal independent set and has a maximum dual value of $1$. If the number of $2$-hop neighbors for a variable is large, the variable can appear in many independent sets, so this variable tends to have a small dual solution value because all constraints corresponding to these independent sets in Eq.~\eqref{eq:dp-2} must be satisfied.


\subsection{Adaptive Stabilization}
\label{subsec:apm}

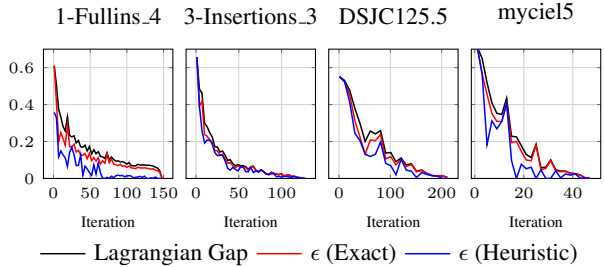
\begin{figure}[t!]
\centering
	\begin{tikzpicture}
	\begin{groupplot}[group style = {group size = 5 by 1, horizontal sep = 5pt}, height=0.4\columnwidth, width=0.4\columnwidth, grid style={line width=.1pt, draw=gray!10},major grid style={line width=.2pt,draw=gray!30}, xlabel = \tiny Iteration, ticklabel style = {font=\tiny}, ymin=0, ymax=0.7, xmajorgrids=true, ymajorgrids=true,  major tick length=0.05cm, minor tick length=0.0cm, legend style={font=\tiny, column sep = 1pt, legend columns = 3,draw=none}]
 
   \nextgroupplot[%
    legend to name=group1,
    title=\small 1-Fullins\_4,
    ]

	\addplot[color=black, line width=0.2mm] table [x index={0}, y index={1}, col sep=space, each nth point=3] {img/ascg/1-FullIns_4.dat};\addlegendentry{\small Lagrangian Gap}
	\addplot[color=red, line width=0.2mm] table [x index={0}, y index={2}, col sep=space, each nth point=3] {img/ascg/1-FullIns_4.dat};\addlegendentry{\small $\epsilon$ (Exact)}
 \addplot[color=blue, line width=0.2mm] table [x index={0}, y index={3}, col sep=space, each nth point=3] {img/ascg/1-FullIns_4.dat};\addlegendentry{\small $\epsilon$ (Heuristic)}

    \nextgroupplot[%
    title=\small 3-Insertions\_3,
    yticklabel=\empty
    ]
	\addplot[color=black, line width=0.2mm] table [x index={0}, y index={1}, col sep=space, each nth point=3] {img/ascg/3-Insertions_3.dat};
	\addplot[color=red, line width=0.2mm] table [x index={0}, y index={2}, col sep=space, each nth point=3] {img/ascg/3-Insertions_3.dat};
 \addplot[color=blue, line width=0.2mm] table [x index={0}, y index={3}, col sep=space, each nth point=3] {img/ascg/3-Insertions_3.dat};

     \nextgroupplot[%
    title=\small DSJC125.5,
    yticklabel=\empty
    ]
	\addplot[color=black, line width=0.2mm] table [x index={0}, y index={1}, col sep=space, each nth point=10] {img/ascg/DSJC125.5.dat};
	\addplot[color=red, line width=0.2mm] table [x index={0}, y index={2}, col sep=space, each nth point=10] {img/ascg/DSJC125.5.dat};
 \addplot[color=blue, line width=0.2mm] table [x index={0}, y index={3}, col sep=space, each nth point=10] {img/ascg/DSJC125.5.dat};

     \nextgroupplot[%
    title=\small myciel5,
    yticklabel=\empty
    ]
	\addplot[color=black, line width=0.2mm, each nth point=2] table [x index={0}, y index={1}, col sep=space] {img/ascg/myciel5.dat};
	\addplot[color=red, line width=0.2mm, each nth point=2] table [x index={0}, y index={2}, col sep=space] {img/ascg/myciel5.dat};
 \addplot[color=blue, line width=0.2mm, each nth point=2] table [x index={0}, y index={3}, col sep=space] {img/ascg/myciel5.dat};
 
    \end{groupplot} 
    \node at (group c1r1.south) [anchor=south, yshift=-1.4cm, xshift=2.7cm] {\pgfplotslegendfromname{group1}}; 
	\end{tikzpicture}
\caption{The trend of penalty value $\epsilon$ computed by exact or heuristic pricing in different iterations of ASCG, both bounded by the Lagrangian Gap on several Graph Coloring Benchmarks.}
\label{fig:t}
\end{figure}

We introduce an adaptive stabilization method to optimize the use of ML prediction. This method adjusts the penalty value $\epsilon$ according to the progress of CG and ensures that the value remains below a Lagrangian gap. As the quality of dual values improves, the Lagrangian gaps become smaller and so are the values of $\epsilon$. Figure~\ref{fig:t} illustrates this dynamic when ASCG is applied to several graphs, which will be discussed in more detail after introducing our method.

Acquiring such dynamics of the penalty parameter is motivated by observations that dual values in the initial stage of CG can be significantly inaccurate (Figure~\ref{fig:cg}), which justifies a stronger influence (i.e., a greater penalty) of the ML predictions to expedite the improvement of dual values. As CG advances, the accuracy of dual values improves, eventually exceeding that of the ML prediction. Consequently, it becomes necessary to decrease the penalty parameter to facilitate the convergence of CG in its final stage. 

Specifically, we propose to set $\epsilon$ as follows:
\begin{equation}
\label{eq:epsilon}
    \epsilon = 
    \begin{cases}
        \frac{c^*_{\epsilon}}{c^*_{\epsilon} - 1}, \;\;\;\;\;\;\;\;& \text{if } c^*_{\epsilon} < 0, \\
        0, & \text{if } c^*_{\epsilon} \geq 0.
    \end{cases}
\end{equation}
Here, $c^*_{\epsilon}$ denotes the optimal objective value (i.e., the minimum reduced cost) of the pricing problem, which is constructed using dual values $\bm{\pi}^{\epsilon}$ to the G-RDP in the current CG iteration. The updated value $\epsilon$ will be used to define the G-RDP in the subsequent iteration of CG. 

\paragraph{Exact Pricing.} We show that the value $\epsilon$ produced by Eq.~\eqref{eq:epsilon} is always upper bounded by a Lagrangian gap in Theorem~\ref{theorem:bound}. The Lagrangian gap measures the quality of the current dual values and reveals the progress of CG. Our definition of the Lagrangian gap is based on the notion of the Lagrangian bound, which provides a valid lower bound to the optimal objective value of the dual problem~\citep{lubbecke2005selected}. Let $z$ denote the objective function of the dual problem, Eq.~\eqref{eq:dualobj}. Given dual values to the G-RDP ($\bm{\pi}^{\epsilon}$), the Lagrangian dual bound can be computed as:
\begin{equation}
    L(\bm{\pi}^{\epsilon}) = \frac{z(\bm{\pi}^{\epsilon})}{(1-c^*_{\epsilon})}.
\end{equation}

\begin{theorem}\label{theorem:bound}
    Define the Lagrangian gap as $l = 1 - L(\bm{\pi}^{\epsilon}) / z(\bm{\pi})$, where $\bm{\pi}$ denotes the optimal dual values associated with the standard RDP. The penalty value $\epsilon$ is in the range $0 \le \epsilon \le l$.
\end{theorem}

\paragraph{Extension to Heuristic Pricing.} We have established results for adaptive stabilization. However, we have assumed that the minimum reduced cost is available, that is, the pricing problem can be solved to optimality. This is not always the case as pricing problems are typically NP-hard. It is common in practice to find heuristic solutions to pricing problems, and an exact pricing solver is only necessary to be used once a heuristic pricing method cannot find any solution with a negative reduced cost. Although using a heuristic pricing method typically increases the number of iterations required for CG to converge, the overall computational time may be reduced because of more efficient solving of the pricing problems.

Therefore, we investigate the effect of heuristic pricing on the proposed adaptive stabilization method. Given the dual values of a G-RDP and its associated pricing problem, we show that the penalty computed using a heuristic pricing method is upper bounded by that with an exact pricing method in Lemma~\ref{lemma:2}. Using Theorem~\ref{theorem:bound}, it is easy to see that the adaptive stabilization method with heuristic pricing is also upper bounded by the Lagrangian gap (Corollary~\ref{coro:1}).

\begin{lemma}
\label{lemma:2}
The penalty value $\epsilon^{\ddagger}$ produced by ASCG under the heuristic pricing setting is upper bounded by the penalty value $\epsilon$ produced by ASCG under the exact pricing, $\epsilon^{\ddagger} \le \epsilon$.

\end{lemma}

\begin{corollary}
\label{coro:1}
The penalty value, $\epsilon^{\ddagger}$, calculated using the adaptive stabilization method with a heuristic method, is upper bounded by a Lagrangian gap, $l$. 
\end{corollary}

Figure~\ref{fig:t} shows the dynamics of penalty values computed with the adaptive stabilization method using an exact or a heuristic pricing method. First, it can be observed that the heuristic penalty curve is always below the exact penalty curve, which is, in turn, below the Lagrangian gap curve. Second, the exact penalty curve closely resembles the Lagrangian gap curve, showing that the exact penalty curve indeed reflects the current progress of CG to a great extent. Compared to that, the heuristic penalty is less reliable, as indicated by the heavier fluctuation of the corresponding curve in some cases. This may lead to an early degeneration of our ASCG-ML to the standard CG method. Also, note that the Lagrangian gap does not decrease monotonically, and computing it requires a dual solution to the standard RDP. This information is not available in ASCG and has to be computed with a substantial additional cost. 
\paragraph{Convergence Property.} Finally, we show in Lemma~\ref{lem:convergence} that the parameter $\epsilon$ in our ASCG method converges to zero after a finite number of iterations. At this stage, our ASCG method becomes equivalent to the standard CG and, therefore, inherits the convergence guarantee of CG.
\begin{lemma}\label{lem:convergence}
     The adaptive stabilization method ensures that $\epsilon$ reduces to $0$ in a finite number of iterations. 
\end{lemma}

\section{Empirical Studies}

\subsection{Setup}

We consider two graph benchmarks: the Matilda library~\citep{smith2014towards} and the Graph Coloring Benchmarks\footnote{\url{https://sites.google.com/site/graphcoloring}} (GCB). We will use Matilda for training and use both benchmarks for testing.  Matilda consists of a large set of $8278$ graphs with $100$ vertices, providing sufficient training data. These graphs are generated following a systematic procedure: an instance-space analysis is carried out using test graphs collected from past studies for GCP and DIMACS graphs; then, an evolutionary algorithm is used to evolve random graphs with different characteristics. 

Our training data are collected from $1258$ Matilda graphs that cannot be solved by CG in $100$ iterations. We randomly select $800$ graphs for training $200$ graphs for validation, and the remaining graphs for testing. When constructing the training data, we set the prediction target for a variable to the average of its optimal dual values when the LP associated with a graph has multiple optimal dual solutions. Our statistical features are computed on a set of $5n$ randomly sampled MISs, where $n$ is the number of vertices in a graph. We employ commonly used training procedures and hyperparameters to train a $3$-layer FFNN and a $20$-layer GCN. In particular, we note that an early stopping criterion and an $L_2$ regularization are used for training to improve the generalization of ML methods. 

The following methods are compared in this study: 

\textit{ASCG-ML} (Ours). We implement our methods according to Algorithm~\ref{alg:ada}. G-RDP is solved by the LP solver in Gurobi~\citep{gurobi2018gurobi} and the pricing problem is solved by a specialized exact solver, TSM~\citep{jiang2018two}. We refer to ASCG with FFNN and GCN as ASCG-FFNN and ASCG-GCN, respectively. 

\textit{SCG-ML}~\citep{kraul2023machine}: This method can be viewed as a special case of the ASCG-ML method with a constant penalty value $\epsilon = 0.1$, selected from $\{0.01,0.1,1,10\}$. As SCG-ML was originally proposed for the cutting stock problem, we equip it with our ML models. We denote SCG-ML using FFNN and GCN as SCG-FFNN and SCG-GCN. 

\textit{SCG/ASCG-Deg}: This method uses the degree attribute to estimate the dual values for SCG or ASCG. We include this baseline to examine the need to develop advanced ML models for dual solution prediction. Note that the degree attribute is used as an estimate of good dual values in a previous study~\citep{briant2008comparison}. 

\textit{SCG}~\citep{du1999stabilized}: This method aims to mitigate the oscillation between successive dual iterates. It performs ASCG with the ML prediction $\hat{\bm{y}}$ replaced by the dual values in the previous iteration, and uses a constant penalty value $\epsilon = 1$, selected from $\{0.01,0.1,1,10\}$.

\textit{Classic CG}~\citep{mehrotra1996column}: The classic CG is a special case of our method with $\epsilon=0$. 

For all compared methods, we provide them an effective warm start with a widely used strategy~\citep{lubbecke2005selected}. Specifically, we use columns from a high-quality integer solution to initialize G-RDP. The heuristic solution is produced by a greedy heuristic and a local search method~\citep{galinier2006survey}. For all stabilization methods, we set the penalty parameter directly to $0$ when its value is below $0.01$, to avoid numerical instability and ensure CG convergence. For all SCG approaches, their penalty values are reduced by half when the pricing problem cannot produce any column with a negative reduced cost.

We note that the complete experiment setup can be found in the appendix. These include the graph benchmark statistics (Appendix~\ref{appendix:gstats}), empirical studies of the labeling approach (Appendix~\ref{appendix:label}), training procedure and layer tuning for ML models (Appendix~\ref{appendix:mlsetup}), and tuning of the penalty values for the compared stabilization methods (Appendix~\ref{appendix:epstune}). 

For testing, we conduct $10$ independent runs for each method on each graph with a random seed in $[1,2,\cdots,10]$. Random seeds can have an impact on the warm-start method and the heuristic pricing method. We measure the convergence rate of a method by its iteration number to reach optimality and report the computational time in seconds. All results are averaged using the geometric mean. Our experiments are conducted on a cluster with $24$ CPUs (Intel(R) Xeon(R) Gold 6238R CPU @ 2.20GHz). Our code is available at \url{https://github.com/yunzhuangshen/ML-based-Adaptive-Stabilization}.

\subsection{Results on Matilda}

\begin{table}[t!]
    \centering
    \caption{Results on Matilda.} 
    \label{tab:cg-matilda}
        \resizebox{0.98\columnwidth}{!}{
        \begin{tabular}{@{}lccccccc@{}}
      \toprule
      & \# Iter. & \% Red. & \# Win & Time & \% Red. & \# Win & P Time\\
      \cmidrule(lr){1-1} \cmidrule(lr){2-4}\cmidrule(lr){5-7}\cmidrule(lr){8-8}
      CG & 341.1 & N/A & 0 & 1921 & N/A & 1 & 1.0\\
      \cmidrule(lr){1-1} \cmidrule(lr){2-4}\cmidrule(lr){5-7}\cmidrule(lr){8-8}
      SCG & 289.3 & 15\% & 0 &  1691 & 12\% & 0 & 1.1\\        
              
      SCG-Deg & 318.1 & 6\% & 0 & 2748 & -43\% & 0 & 1.11\\
      SCG-FFNN & 194.3 & 43\% & 0 & 1450 & 24\% & 1 & 1.32\\
      SCG-GCN & 195.6 & 43\% & 0 & 1463 & 24\% & 3 & 1.34\\            

    \cmidrule(lr){1-1} \cmidrule(lr){2-4}\cmidrule(lr){5-7}\cmidrule(lr){8-8}

      ASCG-Deg & 309.5 & 9\% & 0 & 2467 & -28\% & 0 & 1.06\\
      ASCG-FFNN & \textbf{163.0} & \textbf{52\%} &  \textbf{165} & \textbf{1103} & \textbf{43\%} & \textbf{149} & 1.27\\
      ASCG-GCN & 166.5 & 51\% & 93 & 1130 & 41\% & 104 & 1.27\\   
        \bottomrule
    \end{tabular}
    }
\end{table}

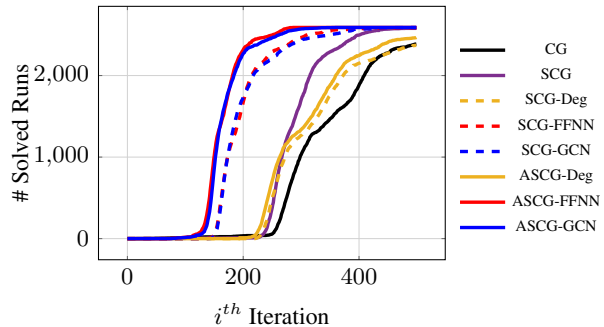
\begin{figure}[t!]
    \centering
	\begin{tikzpicture}
	\begin{axis}[height=0.6\columnwidth, width=0.75\columnwidth, grid style={line width=.1pt, draw=gray!10},major grid style={line width=.2pt,draw=gray!30}, xmajorgrids=true, ymajorgrids=true,  major tick length=0.05cm, minor tick length=0.0cm, ylabel = \small \# Solved Runs, y label style={at={(axis description cs:-0.175,.5)},anchor=south}, xlabel = \small $i^{th}$ Iteration, 
legend style={at={(1.25,0.9)}, anchor=north, font=\tiny, draw=none}
 , ticklabel style = {font=\small}]
	\addplot[color=black, line width=0.45mm] table [x=x, y=CG, col sep=comma]
	{img/Matilda-iter-solved.csv};\addlegendentry{CG}
	\addplot[color=purple, line width=0.45mm] table [x=x, y=SCG, col sep=comma] {img/Matilda-iter-solved.csv};\addlegendentry{SCG}
 
        \addplot[color=orange, line width=0.45mm, dashed] table [x=x, y=SCG-Deg, col sep=comma]
	{img/Matilda-iter-solved.csv};\addlegendentry{SCG-Deg}
        \addplot[color=red, line width=0.45mm, dashed] table [x=x, y=SCG-FFNN, col sep=comma]
	{img/Matilda-iter-solved.csv};\addlegendentry{SCG-FFNN}
        \addplot[color=blue, line width=0.45mm, dashed] table [x=x, y=SCG-GCN, col sep=comma]
	{img/Matilda-iter-solved.csv};\addlegendentry{SCG-GCN}
        \addplot[color=orange, line width=0.45mm] table [x=x, y=ASCG-Deg, col sep=comma]
	{img/Matilda-iter-solved.csv};\addlegendentry{ASCG-Deg}
        \addplot[color=red, line width=0.45mm] table [x=x, y=ASCG-FFNN, col sep=comma]
	{img/Matilda-iter-solved.csv};\addlegendentry{ASCG-FFNN}
        \addplot[color=blue, line width=0.45mm] table [x=x, y=ASCG-GCN, col sep=comma]
	{img/Matilda-iter-solved.csv};\addlegendentry{ASCG-GCN}
 
    \end{axis} 
	\end{tikzpicture}
    \caption{The number of solved runs within an iteration limit on the Matilda graphs.}
    \label{fig:matilda}
\end{figure}

Table~\ref{tab:cg-matilda} shows the results for CG. Here, `\# Iter' denotes the average number of iterations. `\# Time' denotes the average computation time for a seeded run on all $258$ Matilda graphs. `\% Red.' measures the percentage of iterations (or time) reduced by a method compared to the baseline CG. `\# Win' shows the number of graphs (among 258 graphs) that a method solved using the smallest iteration number (or the shortest time). The last column shows the normalized time used for pricing per iteration. 

Overall, the proposed ASCG-FFNN/GCN achieves the most significant reduction in both iteration number and computational time, and they are the winning methods on the majority of graphs according to either evaluation metric. In contrast, ASCG with the basic degree attribute is inefficient, evidencing the need to develop advanced ML techniques that can model the variable interrelationships for accurate predictions. The proposed adaptive stabilization method also shows efficacy. For example, compared to SCG-FFNN, ASCG-FFNN achieves a 9\% further reduction in the number of iterations, while it achieves a more significant 19\% further reduction in computational~time. Consistent with the above analysis, Figure~\ref{fig:matilda} shows that our ASCG-FFNN/GCN method solves more problem instances within a certain iteration cutoff point.

An interesting observation from Table~\ref{tab:cg-matilda} is that the application of stabilization techniques can lead to a more significant reduction in the number of iterations compared to the computational time. The reason lies in the increased difficulty in solving pricing problems, as indicated by the increase in pricing time (the last column of Table~\ref{tab:cg-matilda}). Specifically, stabilization requires solving G-RDP, resulting in more nonzero dual values compared with solving RDP; since the dual values are used as objective coefficients in the pricing problem, the size of the pricing problem increases. As a result, computational time can only be reduced when the reduction in the number of iterations outweighs the computational overhead from solving more difficult pricing problems.

\subsection{Results on Graph Coloring Benchmarks}

In this part, we investigate whether ASCG-ML, trained using Matilda graphs, can be practically used to solve GCB graphs. GCB contains several classes of graphs, each of which has distinct patterns. Moreover, graphs of the same class can vary significantly in size and density. We use $7$ classes of graphs, each with a varying number of graphs: `DSJC' ($12$ graphs), `FullIns' ($14$ graphs), `Insertions' ($11$ graphs), `flat' ($6$ graphs), `le' ($12$ graphs), `myciel' ($5$ graphs), and `queen' ($13$ graphs).  

We evaluate the compared methods with an exact/heuristic pricing method with a cutoff limit of $1000$ seconds. The heuristic pricing setting is included to examine its potential for addressing the issue of more difficult pricing problems in stabilization. We employ a specialized local search method, LSCC~\citep{wang2016two}, as the heuristic pricing method. 

\begin{table}[tb!]
    \centering
    \caption{Results on Graph Coloring Benchmarks.}
    \label{tab:gcb}
    \resizebox{0.96\columnwidth}{!}{\begin{tabular}{lrrrr|rr}
    \toprule
    \multirow{2}{*}{Method} & \multicolumn{4}{c}{Exact Pricing} & \multicolumn{2}{c}{Heuristic Pricing} \\
    \cmidrule(lr){2-5}\cmidrule(lr){6-7}
        &  \# Iter. & \% Red. & Time &  \% Red. & \# Iter. & Time \\
  \cmidrule(lr){1-1} \cmidrule(lr){2-5}\cmidrule(lr){6-7}
           CG &  155.1 &  N/A & 17.8 & N/A &  215.9 & 18.0 \\
  \cmidrule(lr){1-1} \cmidrule(lr){2-5}\cmidrule(lr){6-7}
          SCG &  173.8 &  -12.1\% & 26.1 & -46.6\% & 245.4 &  18.3 \\
     SCG-FFNN &  112.4 &  27.5\% & 18.6 & -4.5\% & 176.5 & 16.5 \\
      SCG-GCN &  143.9 &  7.2\% & 22.6 & -27.0\% & 192.9 & 17.6 \\
  \cmidrule(lr){1-1} \cmidrule(lr){2-5}\cmidrule(lr){6-7}
    ASCG-FFNN &   \textbf{86.6} &  \textbf{44.2}\% & \textbf{16.1} & \textbf{9.6}\% & \textbf{157.4} & \textbf{15.6} \\
     ASCG-GCN &  122.9 &  20.8\% & 18.9 & -6.2\% &  184.5 & 17.2 \\
    \bottomrule
    \end{tabular}}
\end{table}

Table~\ref{tab:gcb} shows the results averaged over all GCB graphs. Specifically, results in `\# Iteration' are averaged over $30$ graphs which can be solved by all compared methods within the cutoff time; results in `Time' are averaged over $42$ graphs, where additional $12$ graphs are those that can be solved by at least one method. When a method does not solve a graph within the $1000$-second cutoff time, this cutoff time is used to produce average results. Note that numerical results on individual Graph Coloring Benchmarks can be found in Appendix~\ref{appendix:ret}.

Overall, we can observe that ASCG-FFNN achieves the best performance in both exact and heuristic pricing setups. First, we focus on the exact pricing setup. The results show that most stabilization approaches can reduce the number of iterations and accelerate the convergence of CG. However, ASCG-FFNN is the only one that outperforms CG in solution time because of the trade-off discussed before. Furthermore, ASCG-FFNN converges with the smallest number of iterations on each class of test graphs, as seen in Table~\ref{tab:gcb_iter}.

Similar observations can be obtained in the heuristic pricing setting. Comparing the results of this setting with the exact pricing setting, all methods require more iterations to converge because of the reduced quality of generated columns, but the computational time may be reduced due to the improved efficiency of pricing. In particular, ML-based stabilization methods achieve a further reduction in computational time, showing that the issue of more difficult pricing problems caused by stabilization can be addressed with a more efficient pricing method.

\begin{table}[tb!]
    \caption{Average iteration number for each class of test graphs. The number of solved graphs in each class is shown in parentheses. -F stands for -FFNN and -G stands for -GCN.}
    \label{tab:gcb_iter}
    \centering
    \resizebox{0.98\columnwidth}{!}{\begin{tabular}{@{}lrrrrrr@{}}
        \toprule
        Graph Class & CG & SCG & SCG-F & SCG-G & ASCG-F & ASCG-G \\
        \midrule
        DSJC (6) & 560.4 & 507.1 & 379.9 & 402.4 & \textbf{349.5} & 368.7 \\
        FullIns (6) & 63.2 & 105.2 & 65.1 & 86.7 & \textbf{31.9} & 58.6 \\
        Insertions (4) & 358.3 & 176.4 & 150.5 & 281.5 & \textbf{133.5} & 269.1 \\
        flat (1) & 1240.9 & 816.4 & 832.8 & 975.5 & \textbf{770.1} & 954.2 \\
        myciel (4) & 71.8 & 60.9 & 48.2 & 55.4 & \textbf{38.2} & 49.9 \\
        queen (9) & 91.5 & 158.0 & 73.0 & 92.6 & \textbf{61.0} & 80.5 \\
        \bottomrule
    \end{tabular}}
\end{table}

\subsection{Discussions \& Analysis on ML Generalization}
\label{subsec:ml}

In this part, we present an analysis by contrasting the performances of ML models for dual solution prediction and the results of ASCG equipped with these methods, to shed light on future advancements of ML techniques. Recall the key difference between the two ML approaches: The GCN approach takes as input a graph encoding all possible interrelationships of dual variables, whereas the FFNN approach utilizes a set of statistical features corresponding to several predefined relations among variables.

\begin{table}[t!]
    \centering
    \caption{Performances of the ML models for dual solution prediction and CG on two graph benchmarks. The ML models are trained using Matilda graphs.}
    \label{tab:ml}
    \resizebox{0.95\columnwidth}{!}{
    \begin{tabular}{lcccccc}
    
    \toprule
             & \multicolumn{3}{c}{Mean Squared Error (1E-2)} & \multicolumn{3}{c}{\# Iteration of CG} \\
        \cmidrule(lr){2-4}\cmidrule(lr){5-7}
             &   FFNN & GCN & Diff. & ASCG-FFNN & ASCG-GCN & Diff. \\
    \cmidrule(lr){1-1} \cmidrule(lr){2-4}\cmidrule(lr){5-7}
     Matilda &   0.33 & 0.30 & -9\%  & 163.0 & 166.5 &  - 2\% \\
     GCB     &   3.8  & 5.9 &  55\% & 86.6 & 122.9 & -34 \% \\
     \bottomrule
    \end{tabular}}
\end{table}

Table~\ref{tab:ml} shows the mean squared error (MSE) for the dual solution prediction task and the iteration number for ASCG equipping an ML model. The MSE is calculated for Matilda using $259$ test graphs and for GCB using $33$ graphs with optimal dual solutions. First, we observe that GCN obtains a smaller test error than FFNN on Matilda, indicating its ability to capture more nuanced interrelationships of variables compared to the predefined statistical features used by FFNN. However, GCN obtains a much higher test error than FFNN in GCB, indicating that the complex variable interrelationships learned by GCN in Matilda are less applicable to GCB. Consequently, ASCG-FFNN outperforms ASCG-GCN significantly on GCB. It is worthnoting that ASCG-GCN receives a marginally worse performance than ASCG-FFNN on Matilda, showing a slight mismatch between the ML prediction accuracy and CG's performance. Specifically, ML models are trained by minimizing the MSE across all training examples associated with dual variables in different problem instances, while the performance of ASCG depends on the quality of the predicted dual solution as a whole.

Nevertheless, we are optimistic about the potential of graph neural networks (GNNs) to expand the applicability of ASCG-ML to a more diverse set of problems. Specifically, GNNs, such as the bipartite graph neural network~\citep{gasse2019exact}, can input a generic description of an optimization problem -- its mixed-integer programming formulation. This allows it to automatically process the rich information in the constraint matrix of a pricing problem and to make predictions for general CG applications. The effectiveness of GNNs could be further enhanced by more closely aligning the training objective with the optimization goal, developing more advanced generalization or regularization strategies such as bilevel optimization~\citep{geisler2021generalization, wang2021bi}, and bridging the gap between training and test data distributions using instance space analysis~\citep{smith2023instance}, all of which merit further research.

\section{Conclusion}

This paper explores machine learning (ML) to enhance Column Generation (CG).  ML is used to make accurate predictions of the optimal dual solution, which is then used to accelerate the convergence of dual values in CG. In addition, we introduce an adaptive stabilization method to make effective use of ML predictions. Our method is termed adaptive stabilized CG based on ML, ASCG-ML.

We show several general properties of ASCG-ML and empirically validate its effectiveness in the graph coloring problem. The findings reveal that ASCG-ML substantially surpasses traditional methods in handling synthetic graphs and exhibits impressive generalization capabilities across a wider range of graphs with distinct patterns. 

Looking ahead, our future work includes investigating ASCG-ML for a broader range of problems and integrating it with more sophisticated elements, such as other stabilization techniques and heuristic pricing methods that can efficiently generate multiple columns.


\section*{Impact Statements}
This paper presents advances in the integration of optimization with machine learning. After careful consideration, we conclude that our work is not expected to cause immediate negative ethical and societal consequences that would require specific emphasis in this context.

\bibliographystyle{icml2024}
\bibliography{lib}

\begin{thebibliography}{55}
\providecommand{\natexlab}[1]{#1}
\providecommand{\url}[1]{\texttt{#1}}
\expandafter\ifx\csname urlstyle\endcsname\relax
  \providecommand{\doi}[1]{doi: #1}\else
  \providecommand{\doi}{doi: \begingroup \urlstyle{rm}\Url}\fi

\bibitem[Agarwal et~al.(1989)Agarwal, Mathur, and Salkin]{agarwal1989set}
Agarwal, Y., Mathur, K., and Salkin, H.~M.
\newblock A set-partitioning-based exact algorithm for the vehicle routing
  problem.
\newblock \emph{Networks}, 19\penalty0 (7):\penalty0 731--749, 1989.

\bibitem[Amor et~al.(2009)Amor, Desrosiers, and Frangioni]{amor2009choice}
Amor, H. M.~B., Desrosiers, J., and Frangioni, A.
\newblock On the choice of explicit stabilizing terms in column generation.
\newblock \emph{Discrete Applied Mathematics}, 157\penalty0 (6):\penalty0
  1167--1184, 2009.

\bibitem[Babaki et~al.(2021)Babaki, Jena, and Charlin]{babaki2021neural}
Babaki, B., Jena, S.~D., and Charlin, L.
\newblock Neural column generation for capacitated vehicle routing.
\newblock In \emph{AAAI-22 Workshop on Machine Learning for Operations Research
  (ML4OR)}, 2021.

\bibitem[Barnhart et~al.(1998)Barnhart, Johnson, Nemhauser, Savelsbergh, and
  Vance]{barnhart1998branch}
Barnhart, C., Johnson, E.~L., Nemhauser, G.~L., Savelsbergh, M. W.~P., and
  Vance, P.~H.
\newblock Branch-and-price: Column generation for solving huge integer
  programs.
\newblock \emph{Oper. Res.}, 46\penalty0 (3):\penalty0 316--329, 1998.

\bibitem[Ben~Amor et~al.(2006)Ben~Amor, Desrosiers, and Val{\'e}rio~de
  Carvalho]{ben2006dual}
Ben~Amor, H., Desrosiers, J., and Val{\'e}rio~de Carvalho, J.~M.
\newblock Dual-optimal inequalities for stabilized column generation.
\newblock \emph{Operations Research}, 54\penalty0 (3):\penalty0 454--463, 2006.

\bibitem[Bengio et~al.(2021)Bengio, Lodi, and Prouvost]{bengio2020machine}
Bengio, Y., Lodi, A., and Prouvost, A.
\newblock Machine learning for combinatorial optimization: {A} methodological
  tour d'horizon.
\newblock \emph{Eur. J. Oper. Res.}, 290\penalty0 (2):\penalty0 405--421, 2021.

\bibitem[Briant et~al.(2008)Briant, Lemar{\'e}chal, Meurdesoif, Michel, Perrot,
  and Vanderbeck]{briant2008comparison}
Briant, O., Lemar{\'e}chal, C., Meurdesoif, P., Michel, S., Perrot, N., and
  Vanderbeck, F.
\newblock Comparison of bundle and classical column generation.
\newblock \emph{Mathematical programming}, 113\penalty0 (2):\penalty0 299--344,
  2008.

\bibitem[Chi et~al.(2022)Chi, Aboussalah, Khalil, Wang, and
  Sherkat-Masoumi]{chi2022deep}
Chi, C., Aboussalah, A., Khalil, E., Wang, J., and Sherkat-Masoumi, Z.
\newblock A deep reinforcement learning framework for column generation.
\newblock \emph{Advances in Neural Information Processing Systems},
  35:\penalty0 9633--9644, 2022.

\bibitem[Dakin(1965)]{dakin1965tree}
Dakin, R.~J.
\newblock A tree-search algorithm for mixed integer programming problems.
\newblock \emph{The computer journal}, 8\penalty0 (3):\penalty0 250--255, 1965.

\bibitem[Dantzig(2016)]{dantzig2016linear}
Dantzig, G.
\newblock \emph{Linear programming and extensions}.
\newblock Princeton university press, 2016.

\bibitem[Deza \& Khalil(2023)Deza and Khalil]{deza2023machine}
Deza, A. and Khalil, E.~B.
\newblock Machine learning for cutting planes in integer programming: A survey.
\newblock \emph{arXiv preprint arXiv:2302.09166}, 2023.

\bibitem[Ding et~al.(2020)Ding, Zhang, Shen, Li, Wang, Xu, and
  Song]{ding2020accelerating}
Ding, J., Zhang, C., Shen, L., Li, S., Wang, B., Xu, Y., and Song, L.
\newblock Accelerating primal solution findings for mixed integer programs
  based on solution prediction.
\newblock In \emph{The Thirty-Fourth {AAAI} Conference on Artificial
  Intelligence, New York, NY, USA, February 7-12, 2020}, pp.\  1452--1459.
  {AAAI} Press, 2020.

\bibitem[Du~Merle et~al.(1999)Du~Merle, Villeneuve, Desrosiers, and
  Hansen]{du1999stabilized}
Du~Merle, O., Villeneuve, D., Desrosiers, J., and Hansen, P.
\newblock Stabilized column generation.
\newblock \emph{Discrete Mathematics}, 194\penalty0 (1-3):\penalty0 229--237,
  1999.

\bibitem[Formanowicz \& Tana{\'s}(2012)Formanowicz and
  Tana{\'s}]{formanowicz2012survey}
Formanowicz, P. and Tana{\'s}, K.
\newblock A survey of graph coloring-its types, methods and applications.
\newblock \emph{Foundations of Computing and Decision Sciences}, 37\penalty0
  (3):\penalty0 223--238, 2012.

\bibitem[Galinier \& Hertz(2006)Galinier and Hertz]{galinier2006survey}
Galinier, P. and Hertz, A.
\newblock A survey of local search methods for graph coloring.
\newblock \emph{Computers \& Operations Research}, 33\penalty0 (9):\penalty0
  2547--2562, 2006.

\bibitem[Gasse et~al.(2019)Gasse, Ch{\'{e}}telat, Ferroni, Charlin, and
  Lodi]{gasse2019exact}
Gasse, M., Ch{\'{e}}telat, D., Ferroni, N., Charlin, L., and Lodi, A.
\newblock Exact combinatorial optimization with graph convolutional neural
  networks.
\newblock In \emph{Advances in Neural Information Processing Systems 32,
  December 8-14, 2019, Vancouver, BC, Canada}, pp.\  15554--15566, 2019.

\bibitem[Geisler et~al.(2021)Geisler, Sommer, Schuchardt, Bojchevski, and
  G{\"u}nnemann]{geisler2021generalization}
Geisler, S., Sommer, J., Schuchardt, J., Bojchevski, A., and G{\"u}nnemann, S.
\newblock Generalization of neural combinatorial solvers through the lens of
  adversarial robustness.
\newblock \emph{arXiv preprint arXiv:2110.10942}, 2021.

\bibitem[Gilmore \& Gomory(1961)Gilmore and Gomory]{gilmore1961linear}
Gilmore, P.~C. and Gomory, R.~E.
\newblock A linear programming approach to the cutting-stock problem.
\newblock \emph{Operations research}, 9\penalty0 (6):\penalty0 849--859, 1961.

\bibitem[Goffin et~al.(1992)Goffin, Haurie, and Vial]{goffin1992decomposition}
Goffin, J.-L., Haurie, A., and Vial, J.-P.
\newblock Decomposition and nondifferentiable optimization with the projective
  algorithm.
\newblock \emph{Management science}, 38\penalty0 (2):\penalty0 284--302, 1992.

\bibitem[Gondzio \& Sarkissian(1996)Gondzio and Sarkissian]{gondzio1996column}
Gondzio, J. and Sarkissian, R.
\newblock Column generation with a primal-dual method.
\newblock \emph{Relatorio tecnico, University of Geneva}, 102, 1996.

\bibitem[Gondzio et~al.(2013)Gondzio, Gonz{\'a}lez-Brevis, and
  Munari]{gondzio2013new}
Gondzio, J., Gonz{\'a}lez-Brevis, P., and Munari, P.
\newblock New developments in the primal--dual column generation technique.
\newblock \emph{European Journal of Operational Research}, 224\penalty0
  (1):\penalty0 41--51, 2013.

\bibitem[Gschwind \& Irnich(2017)Gschwind and Irnich]{gschwind2017stabilized}
Gschwind, T. and Irnich, S.
\newblock Stabilized column generation for the temporal knapsack problem using
  dual-optimal inequalities.
\newblock \emph{OR spectrum}, 39:\penalty0 541--556, 2017.

\bibitem[Gurobi~Optimization(2018)]{gurobi2018gurobi}
Gurobi~Optimization, I.
\newblock Gurobi optimizer reference manual.
\newblock \url{https://www.gurobi.com}, 2018.
\newblock Accessed: 2022-03-17.

\bibitem[Huang et~al.(2023)Huang, Ferber, Tian, Dilkina, and
  Steiner]{huang2023searching}
Huang, T., Ferber, A.~M., Tian, Y., Dilkina, B., and Steiner, B.
\newblock Searching large neighborhoods for integer linear programs with
  contrastive learning.
\newblock In \emph{International Conference on Machine Learning}, pp.\
  13869--13890. PMLR, 2023.

\bibitem[Jensen \& Toft(2011)Jensen and Toft]{jensen2011graph}
Jensen, T.~R. and Toft, B.
\newblock \emph{Graph coloring problems}.
\newblock John Wiley \& Sons, 2011.

\bibitem[Jiang et~al.(2018)Jiang, Li, Liu, and Many{\`{a}}]{jiang2018two}
Jiang, H., Li, C., Liu, Y., and Many{\`{a}}, F.
\newblock A two-stage maxsat reasoning approach for the maximum weight clique
  problem.
\newblock In \emph{Proceedings of the Thirty-Second {AAAI} Conference on
  Artificial Intelligence, New Orleans, Louisiana, USA, February 2-7, 2018},
  pp.\  1338--1346. {AAAI} Press, 2018.

\bibitem[Khalil et~al.(2016)Khalil, Bodic, Song, Nemhauser, and
  Dilkina]{khalil2016learning}
Khalil, E.~B., Bodic, P.~L., Song, L., Nemhauser, G.~L., and Dilkina, B.
\newblock Learning to branch in mixed integer programming.
\newblock In \emph{Proceedings of the Thirtieth {AAAI} Conference on Artificial
  Intelligence, February 12-17, 2016, Phoenix, Arizona, {USA}}, pp.\  724--731.
  {AAAI} Press, 2016.

\bibitem[Kingma \& Ba(2015)Kingma and Ba]{kingma2014adam}
Kingma, D.~P. and Ba, J.
\newblock Adam: {A} method for stochastic optimization.
\newblock In \emph{3rd International Conference on Learning Representations,
  {ICLR} 2015, San Diego, CA, USA, May 7-9, 2015, Conference Track
  Proceedings}, 2015.

\bibitem[Kipf \& Welling(2017)Kipf and Welling]{kipf2016semi}
Kipf, T.~N. and Welling, M.
\newblock Semi-supervised classification with graph convolutional networks.
\newblock In \emph{5th International Conference on Learning Representations,
  {ICLR} 2017, Toulon, France, April 24-26, 2017, Conference Track
  Proceedings}. OpenReview.net, 2017.

\bibitem[Kool et~al.(2019)Kool, van Hoof, and Welling]{kool2018attention}
Kool, W., van Hoof, H., and Welling, M.
\newblock Attention, learn to solve routing problems!
\newblock In \emph{7th International Conference on Learning Representations,
  New Orleans, LA, USA, May 6-9, 2019}. OpenReview.net, 2019.

\bibitem[Kotary et~al.(2022)Kotary, Fioretto, and
  Van~Hentenryck]{kotary2022fast}
Kotary, J., Fioretto, F., and Van~Hentenryck, P.
\newblock Fast approximations for job shop scheduling: A lagrangian dual deep
  learning method.
\newblock In \emph{Thirty-Sixth {AAAI} Conference on Artificial Intelligence},
  pp.\  7239--7246, 2022.

\bibitem[Kraul et~al.(2023)Kraul, Seizinger, and Brunner]{kraul2023machine}
Kraul, S., Seizinger, M., and Brunner, J.~O.
\newblock Machine learning--supported prediction of dual variables for the
  cutting stock problem with an application in stabilized column generation.
\newblock \emph{INFORMS Journal on Computing}, 2023.

\bibitem[Lee \& Park(2011)Lee and Park]{lee2011chebyshev}
Lee, C. and Park, S.
\newblock Chebyshev center based column generation.
\newblock \emph{Discrete Applied Mathematics}, 159\penalty0 (18):\penalty0
  2251--2265, 2011.

\bibitem[Li et~al.(2018)Li, Chen, and Koltun]{li2018combinatorial}
Li, Z., Chen, Q., and Koltun, V.
\newblock Combinatorial optimization with graph convolutional networks and
  guided tree search.
\newblock In \emph{Advances in Neural Information Processing Systems 31,
  December 3-8, 2018, Montr{\'{e}}al, Canada}, pp.\  537--546, 2018.

\bibitem[Liu et~al.(2022)Liu, Fischetti, and Lodi]{liu2022learning}
Liu, D., Fischetti, M., and Lodi, A.
\newblock Learning to search in local branching.
\newblock In \emph{Thirty-Sixth {AAAI} Conference on Artificial Intelligence},
  pp.\  3796--3803, 2022.

\bibitem[L{\"{u}}bbecke \& Desrosiers(2005)L{\"{u}}bbecke and
  Desrosiers]{lubbecke2005selected}
L{\"{u}}bbecke, M.~E. and Desrosiers, J.
\newblock Selected topics in column generation.
\newblock \emph{Oper. Res.}, 53\penalty0 (6):\penalty0 1007--1023, 2005.

\bibitem[Marsten et~al.(1975)Marsten, Hogan, and
  Blankenship]{marsten1975boxstep}
Marsten, R.~E., Hogan, W.~W., and Blankenship, J.~W.
\newblock The boxstep method for large-scale optimization.
\newblock \emph{Operations Research}, 23\penalty0 (3):\penalty0 389--405, 1975.

\bibitem[Mehrotra \& Trick(1996)Mehrotra and Trick]{mehrotra1996column}
Mehrotra, A. and Trick, M.~A.
\newblock A column generation approach for graph coloring.
\newblock \emph{{INFORMS} J. Comput.}, 8\penalty0 (4):\penalty0 344--354, 1996.

\bibitem[Morabit et~al.(2021)Morabit, Desaulniers, and Lodi]{Morabit2020mlcs}
Morabit, M., Desaulniers, G., and Lodi, A.
\newblock Machine-learning-based column selection for column generation.
\newblock \emph{Transp. Sci.}, 55\penalty0 (4):\penalty0 815--831, 2021.

\bibitem[Pessoa et~al.(2018)Pessoa, Sadykov, Uchoa, and
  Vanderbeck]{pessoa2018automation}
Pessoa, A., Sadykov, R., Uchoa, E., and Vanderbeck, F.
\newblock Automation and combination of linear-programming based stabilization
  techniques in column generation.
\newblock \emph{INFORMS Journal on Computing}, 30\penalty0 (2):\penalty0
  339--360, 2018.

\bibitem[Quesnel et~al.(2022)Quesnel, Wu, Desaulniers, and
  Soumis]{quesnel2022deep}
Quesnel, F., Wu, A., Desaulniers, G., and Soumis, F.
\newblock Deep-learning-based partial pricing in a branch-and-price algorithm
  for personalized crew rostering.
\newblock \emph{Computers \& Operations Research}, 138:\penalty0 105554, 2022.

\bibitem[Rousseau et~al.(2007)Rousseau, Gendreau, and
  Feillet]{rousseau2007interior}
Rousseau, L.-M., Gendreau, M., and Feillet, D.
\newblock Interior point stabilization for column generation.
\newblock \emph{Operations Research Letters}, 35\penalty0 (5):\penalty0
  660--668, 2007.

\bibitem[Shen et~al.(2021)Shen, Sun, Eberhard, and Li]{shen2021learning}
Shen, Y., Sun, Y., Eberhard, A., and Li, X.
\newblock Learning primal heuristics for mixed integer programs.
\newblock In \emph{2021 International Joint Conference on Neural Networks
  (IJCNN)}, pp.\  1--8. IEEE, 2021.

\bibitem[Shen et~al.(2022)Shen, Sun, Li, Eberhard, and
  Ernst]{shen2022enhancing}
Shen, Y., Sun, Y., Li, X., Eberhard, A., and Ernst, A.
\newblock Enhancing column generation by a machine-learning-based pricing
  heuristic for graph coloring.
\newblock In \emph{Proceedings of the AAAI Conference on Artificial
  Intelligence}, pp.\  9926--9934, 2022.

\bibitem[Smith-Miles \& Mu{\~n}oz(2023)Smith-Miles and
  Mu{\~n}oz]{smith2023instance}
Smith-Miles, K. and Mu{\~n}oz, M.~A.
\newblock Instance space analysis for algorithm testing: Methodology and
  software tools.
\newblock \emph{ACM Computing Surveys}, 55\penalty0 (12):\penalty0 1--31, 2023.

\bibitem[Smith-Miles et~al.(2014)Smith-Miles, Baatar, Wreford, and
  Lewis]{smith2014towards}
Smith-Miles, K., Baatar, D., Wreford, B., and Lewis, R.
\newblock Towards objective measures of algorithm performance across instance
  space.
\newblock \emph{Computers \& Operations Research}, 45:\penalty0 12--24, 2014.

\bibitem[Sun et~al.(2021)Sun, Li, and Ernst]{sun2019using}
Sun, Y., Li, X., and Ernst, A.~T.
\newblock Using statistical measures and machine learning for graph reduction
  to solve maximum weight clique problems.
\newblock \emph{{IEEE} Trans. Pattern Anal. Mach. Intell.}, 43\penalty0
  (5):\penalty0 1746--1760, 2021.

\bibitem[Sun et~al.(2022)Sun, Wang, Shen, Li, Ernst, and
  Kirley]{sun2022boosting}
Sun, Y., Wang, S., Shen, Y., Li, X., Ernst, A.~T., and Kirley, M.
\newblock Boosting ant colony optimization via solution prediction and machine
  learning.
\newblock \emph{Computers \& Operations Research}, 143:\penalty0 105769, 2022.

\bibitem[Vanderbeck(2000)]{vanderbeck2000dantzig}
Vanderbeck, F.
\newblock On dantzig-wolfe decomposition in integer programming and ways to
  perform branching in a branch-and-price algorithm.
\newblock \emph{Oper. Res.}, 48\penalty0 (1):\penalty0 111--128, 2000.

\bibitem[Wang et~al.(2021)Wang, Hua, Liu, Zhang, Yan, Qi, Yang, Zhou, and
  Yang]{wang2021bi}
Wang, R., Hua, Z., Liu, G., Zhang, J., Yan, J., Qi, F., Yang, S., Zhou, J., and
  Yang, X.
\newblock A bi-level framework for learning to solve combinatorial optimization
  on graphs.
\newblock \emph{Advances in Neural Information Processing Systems},
  34:\penalty0 21453--21466, 2021.

\bibitem[Wang et~al.(2016)Wang, Cai, and Yin]{wang2016two}
Wang, Y., Cai, S., and Yin, M.
\newblock Two efficient local search algorithms for maximum weight clique
  problem.
\newblock In \emph{Proceedings of the Thirtieth {AAAI} Conference on Artificial
  Intelligence, February 12-17, 2016, Phoenix, Arizona, {USA}}, pp.\  805--811.
  {AAAI} Press, 2016.

\bibitem[Wentges(1997)]{wentges1997weighted}
Wentges, P.
\newblock Weighted dantzig-wolfe decomposition for linear mixed-integer
  programming.
\newblock \emph{International Transactions in Operational Research}, 4\penalty0
  (2):\penalty0 151--162, 1997.

\bibitem[Xin et~al.(2021)Xin, Song, Cao, and Zhang]{xin2021neurolkh}
Xin, L., Song, W., Cao, Z., and Zhang, J.
\newblock Neurolkh: Combining deep learning model with lin-kernighan-helsgaun
  heuristic for solving the traveling salesman problem.
\newblock \emph{Advances in Neural Information Processing Systems},
  34:\penalty0 7472--7483, 2021.

\bibitem[Zhou et~al.(2023)Zhou, Wu, Song, Cao, and Zhang]{zhou2023towards}
Zhou, J., Wu, Y., Song, W., Cao, Z., and Zhang, J.
\newblock Towards omni-generalizable neural methods for vehicle routing
  problems.
\newblock In \emph{International Conference on Machine Learning, {ICML} 2023,
  23-29 July 2023, Honolulu, Hawaii, {USA}}, pp.\  42769--42789, 2023.

\bibitem[Zong et~al.(2022)Zong, Zheng, Li, and Jin]{zong2022mapdp}
Zong, Z., Zheng, M., Li, Y., and Jin, D.
\newblock Mapdp: Cooperative multi-agent reinforcement learning to solve pickup
  and delivery problems.
\newblock In \emph{Thirty-Sixth {AAAI} Conference on Artificial Intelligence},
  pp.\  9980--9988, 2022.

\end{thebibliography}

\newpage
\appendix
\onecolumn

\section{List of Notations}
\label{appendix:notation}
\begin{table}[h!]
    \centering
    \caption{List of Notations.}
    \begin{tabular}{lll}
    \toprule
      \multirow{3}{*}{Graph} & $\mathcal{V}$ & The set of vertices in a graph \\
           & $\mathcal{E}$ & The set of edges in a graph \\
           & $s$ & Set representation of a maximal independent set (MIS)\\
           & $\mathcal{S}$ & A set of maximal independent sets \\
        \midrule
      \multirow{3}{*}{Restricted Dual Problem (RDP)} & $z$ & The objective function  \\
                                                 & $z^*$ & The optimal objective value\\
                                                 & $\bm{\pi}$ & dual variables or dual solution\\
         \midrule
        \multirow{3}{*}{Generalized RMP (G-RDP)} & $z_{\epsilon}$ & The objective function  \\
                                                 & $z^*_{\epsilon}$ & The optimal objective value\\
                                                 & $\bm{\pi}_\epsilon$ & dual variables or dual solution\\                    
        \midrule
      \multirow{2}{*}{Pricing Problem}           & $s^*$ & An optimal solution. \\ 
                                                 & $c^*$ & The minimum reduced cost w.r.t. $\bm{\pi}$\\
                                                 & $c^*_{\epsilon}$ & The minimum reduced cost w.r.t. $\bm{\pi}_{\epsilon}$\\
                                                 & $c^\ddagger$ & The reduced cost of a heuristic solution \\
                                                 & $v$ & The decision variable indicating whether to use a vertex or not \\
        \midrule
      \multirow{4}{*}{Dual Prediction Problem}   & $\mathbb{D}$ & Training data set \\
                                                  & $\bm{f}$ & Feature vector \\
                                                  & $y$ & Prediction target \\
                                                  & $\hat{y}$ & Predicted value \\
          \midrule
        \multirow{3}{*}{Others}                  & $L$ & A Lagrangian bound \\
                                                 & $l$ & A Lagrangian gap \\
                                                 & $\epsilon$ & Penalty value computed with an optimal solution \\
                                                 & $\epsilon_{\ddagger}$ & Penalty value computed with a heuristic solution \\
                                
        \bottomrule                              
    \end{tabular}
\end{table}

\section{The Primal LP Problem for Graph Coloring}
\label{appendix:gcp}

Let $G(\mathcal{V},\mathcal{E})$ be a graph, where $\mathcal{V}$ and $\mathcal{E}$ denote the set of vertices and edges, respectively. Let $x_s$ be a binary variable indicating whether a MIS~$s\in\mathcal{S}$ in the graph is used ($x_s=1$) or not ($x_s=0$). Here, $s$ is a set presentation of a MIS, and $\mathcal{S}$ denotes the set of all MISs in a graph. The Dantizig-Wolfe formulation of the GCP can be defined as:
\begin{align}
    z^* = \min_{\bm{x}} \;& \sum_{s \in \mathcal{S}} x_s, & (\mathcal{P}) \\
        s.t. \;& \sum_{s \in \mathcal{S},i\in s} x_s \geq 1, &  i \in \mathcal{V},  \label{eq:cover}\\
              & x_s \in \{0,1\}, & s\in \mathcal{S}.
\end{align}

\noindent In Eq.~\eqref{eq:cover}, the expression $s\in S, i \in s$ denotes all MISs that contain the vertex $i$. As a graph can contain exponentially many MISs, this LP relaxation can contain many variables (or columns).  The primal LP problem can be obtained by relaxing the integer constraints on $x_s$.

\section{Proofs}
\label{appendix:proof}

\begin{lemma}
    Let $\bm{\pi}^{\epsilon}$ denote the dual solution to the current G-RDP. The minimum reduced cost $c^*_{\epsilon}$ with respect to the dual values $\bm{\pi}^{\epsilon}$ can be positive. 
\end{lemma}
\begin{proof}

     The minimum reduced cost $c^*_{\epsilon}$ can be positive if the solution $\bm{\pi}^{\epsilon}$ to the G-RDP is an interior point in the polyhedron of the original dual problem. In the extreme case where $\epsilon \rightarrow \infty$, $\bm{\pi}^{\epsilon}$ and $\hat{y}$ overlaps. $\bm{\pi}^{\epsilon}$ can be an interior point because ML prediction $\hat{y}$ could be. Another scenario is that the solution of the standard RDP in a CG iteration, $\bm{\pi}$, is on the boundary of the polyhedron of the dual problem and the ML prediction is an interior point of the dual problem. For $\epsilon > 0$, the ML prediction will attract the dual solution to the interior. In a more general case, where the solution $\bm{\pi}$ of the standard RDP in a CG iteration is not on the boundary of the polyhedron of the dual problem, it is still possible that $\bm{\pi}^{\epsilon}$ is an interior point. This is because the ML prediction $\hat{y}$ can be ``moved" further away from $\bm{\pi}$ to attract it to the interior of the dual problem. Essentially, increasing the distance between $\hat{y}$ and $\bm{\pi}$ is equivalent to increasing the value of penalty $\epsilon$.
\end{proof}

\begin{theorem}
    Define the Lagrangian gap as $l = 1 - L(\bm{\pi}^{\epsilon}) / z(\bm{\pi})$, where $\bm{\pi}$ denotes the optimal dual values associated with the standard RDP. The penalty value $\epsilon$ is in the range $0 \le \epsilon \le l$.
\end{theorem}
\begin{proof}
If $c^*_{\epsilon} \ge 0$, $\epsilon  = 0$ according to its definition in Eq.~\eqref{eq:epsilon} and the theorem satisfies. If $c^*_{\epsilon} < 0$, we have
\begin{equation}
\epsilon  = \frac{c^*_{\epsilon}}{c^*_{\epsilon} - 1} > 0. 
\end{equation}
In addition, 
\begin{equation}
     \epsilon  = \frac{c^*_{\epsilon}}{c^*_{\epsilon} - 1}
               = 1- \frac{1}{1-c^*_{\epsilon}} 
               = 1- \frac{\frac{z(\bm{\pi}^{\epsilon})}{1-c^*_{\epsilon}}}{ z(\bm{\pi}^{\epsilon})}
               = 1 - \frac{ L(\bm{\pi_{\epsilon}})}{ z(\bm{\pi}^{\epsilon})}.
\end{equation}
As the optimal solution of the modified RDP, $\bm{\pi}^{\epsilon}$, is a feasible solution to the RDP, the objective value of $\bm{\pi}^{\epsilon}$ in the RDP must be less than or equal to its optimal objective value, $z(\bm{\pi})$, since the RDP is a maximization problem:
\begin{equation}
    z(\bm{\pi}^{\epsilon}) \le z(\bm{\pi}). 
\end{equation}
Hence, 
\begin{equation}
    \epsilon = 1 - \frac{ L(\bm{\pi_{\epsilon}})}{ z(\bm{\pi}^{\epsilon})} \le 1 - \frac{L(\bm{\pi_{\epsilon}})}{z(\bm{\pi})} = l. 
\end{equation}

\end{proof}

\begin{lemma}
The penalty value $\epsilon^{\ddagger}$ produced by ASCG under the heuristic pricing setting is upper bound by the penalty value $\epsilon$ produced by ASCG under the exact pricing setting, $\epsilon^{\ddagger} \le \epsilon$.
\end{lemma}
\begin{proof}
Let $c^{*}$ and $c^{\ddagger}$ denote the reduced costs of an optimal and heuristic solution to the pricing problem, respectively. Since the pricing problem is a minimization problem, we have the following.
\begin{equation}
    c^{*} \le c^{\ddagger}. 
\end{equation}
In the case where both heuristic and exact solutions have negative reduced costs, $c^* \leq c^{\ddagger} < 0$, we have the following. 
\begin{equation}
    \epsilon = \frac{c^{*}}{c^{*}-1} \ge \frac{c^{\ddagger}}{c^{\ddagger}-1} = \epsilon^{\ddagger} > 0. 
\end{equation}
In the case where the reduced costs of both solutions are positive, $0 < c^* \leq c^{\ddagger}$, the penalties are both equal to $0$:
\begin{equation}
    \epsilon = \epsilon^{\ddagger} = 0. 
\end{equation}
In the case where $c^* \le 0 \le c^{\ddagger}$, the adaptive penalty with the optimal solution is nonnegative and that with the heuristic solution is equal to $0$:  
\begin{equation}
    \epsilon \ge 0 = \epsilon^{\ddagger}. 
\end{equation}
\end{proof}

\begin{corollary}
The penalty value, $\epsilon^{\ddagger}$, calculated using the adaptive stabilization method with a heuristic method, is upper bounded by a Lagrangian gap, $l$. 
\end{corollary}
\begin{proof}

Since $\epsilon^{\ddagger} \le \epsilon$ according to Lemma~\ref{lemma:2} and $\epsilon \le l$ according to Theorem~\ref{theorem:bound}, the adaptive penalty $\epsilon^{\ddagger}$ with a heuristic pricing method is upper bound by the Lagrangian gap, $l$.
 \end{proof}

\begin{lemma}
     The adaptive stabilization method ensures that $\epsilon$ reduces to $0$ in a finite number of iterations. 
\end{lemma}
\begin{proof}
    If the minimum reduced cost $c^*_{\epsilon}$ is negative in one CG iteration
    \begin{equation}
        c^*_{\epsilon} =   1 - \sum_{i \in s^*} \pi_i^{\epsilon} < 0, 
    \end{equation}
    the generated column $s^*$ must be a new column that can cut the current optimal solution ($\bm{\pi}^{\epsilon}$) to the modified RDP, as $s^*$ violates the constraint~\eqref{eq:dualcut}. As the number of MISs in a graph is finite, no new columns can be further generated after a finite number of iterations. At this point, the minimum reduced cost is non-negative and, therefore, $\epsilon$ reduces to $0$. 
\end{proof}

\section{Statistical Features}
\label{appendix:mlmodel}
In GCP, a problem instance (i.e., a dual problem) is associated with a graph, a column represents a MIS, and a dual variable corresponds to a vertex. We use two basic graph attributes to describe a variable, the degree and the graph density. The set of statistical features is computed on randomly sampled MISs. Our first statistical feature is the frequency of a vertex~$i$ appearing in the samples $\mathcal{S}$,
\begin{equation}
    f_1(i) = |\mathcal{S}_i|/|\mathcal{S}|,
\end{equation}
where $\mathcal{S}_i$ denotes the set of MISs in $\mathcal{S}$ that contains the vertex~$i$. A vertex with a high frequency indicates that this vertex is likely to be covered, hence it tends to have a small dual solution value.

The second set of statistical features are the maximum, minimum, and average cardinality of MISs containing the vertex $i$ ($\mathcal{S}_i$), defined as,
\begin{equation}
    f_{2-4}(i) = \max_{s\in \mathcal{S}_i}|s|,  \; \min_{ s\in \mathcal{S}_i}|s|, \; \sum_{s\in \mathcal{S}_i}|s|/|\mathcal{S}_i|. 
\end{equation}
A vertex with large cardinality statistics tends to have a small dual solution value, as the sum of the dual solution value for vertices in an MIS should be no greater than $1$, i.e., the reduced cost of a column is non-negative. 

The last set of features is statistics on the average degree of sampled MISs containing that vertex. We denote the average degree of vertices in an MIS $s$ as
\begin{equation}
\overline{\deg}(s) = \sum_{i \in s}\deg(i)/|s|.
\label{eq:avedeg}
\end{equation}
The maximum, minimum, and average statistics of the quantity \eqref{eq:avedeg} across the MISs containing the vertex $i$ ($\mathcal{S}_i$) can be defined respectively as 
\begin{equation}
  \begin{gathered}
    f_{5-7}(i) = \max_{s\in \mathcal{S}_i}\overline{\deg}(s), \; \min_{s\in \mathcal{S}_i}\overline{\deg}(s),  \; \sum_{s\in \mathcal{S}_i}\overline{\deg}(s)/|\mathcal{S}_i|.
  \end{gathered}
\end{equation}
If the average degree of vertices in the MISs containing the vertex $i$ is high, that vertex tends to have a small value. This is because vertices with high degree values tend to have large dual solution values and the sum of the dual solution value for vertices in an MIS should be no greater than $1$. All statistical features are normalized instance-wise, i.e., for a feature the value for a variable is scaled by the maximum and minimum values of variables in the same LP instance.

\section{Details of Experiment Setup}
\label{appendix:expsetup}

\subsection{Benchmark Statistics}
\label{appendix:gstats}

\begin{figure}[tb!]
    \begin{minipage}{0.55\textwidth}
        \centering
              \begin{tikzpicture}
        \begin{axis} [ybar, height=0.6\textwidth,width=0.9\textwidth, xlabel =  Density interval, ylabel =  Count,  axis lines*=left, ytick=\empty,ymin=0,xticklabel style={anchor=base,yshift=-\baselineskip, rotate=30}, xticklabel shift=5pt,
        symbolic x coords={[0-0.1],[0.1-0.2],[0.2-0.3],[0.3-0.4],[0.4-0.5],[0.5-0.6],[0.6-1.0]},nodes near coords={\pgfmathprintnumber\pgfplotspointmeta}]
        \addplot[fill=maroon, draw=maroon ] coordinates {
        ([0-0.1],62) ([0.1-0.2],58) ([0.2-0.3],128) ([0.3-0.4],303) ([0.4-0.5],174) ([0.5-0.6],597) ([0.6-1.0],4)
        }
        \closedcycle;
        \end{axis}
        \end{tikzpicture}
      \captionof{figure}{ Density histogram of the $1326$ Matilda graphs.}
      \label{fig:hist-matilda}
    \end{minipage}
    \hspace{-2em}
    \begin{minipage}{0.4\textwidth}
        \captionof{table}{Statistics of $8$ classes of test graphs from GCB.}
        \label{tab:gcb-stats}
        \centering
\resizebox{\columnwidth}{!}{\begin{tabular}{rrrr}
        \toprule
        Class & \# Instances & \# Nodes & Density  \\
        \midrule
            DSJC & 12 & $[125, 1000]$ & $[0.095, 0.901]$ \\
            FullIns & 14 & $[30, 4146]$ & $[0.009, 0.23]$ \\
            Insertions & 11 & $[37, 1406]$ & $[0.01, 0.108]$ \\
            flat & 6 & $[300, 1000]$ & $[0.477, 0.494]$ \\
            le450 & 12 & $[450, 450]$ & $[0.057, 0.172]$ \\
            myciel & 5 & $[11, 191]$ & $[0.13, 0.364]$ \\
            queen & 13 & $[25, 256]$ & $[0.3879, 0.533]$ \\
        \midrule
         Overall & 73 & $[11,4146]$ & $[0.009, 0.901]$\\
            \bottomrule
        \end{tabular}}
    \end{minipage}
\end{figure}

Graphs from two benchmarks are used in this paper, Matilda and graph coloring benchmarks (GCB). Figure~\ref{fig:hist-matilda} shows the distribution with respect to the density of $1326$ Matilda graphs. Table~\ref{tab:gcb-stats} shows the statistics for the 7 classes of GCB graphs.

\subsection{Empirical studies on Labeling Approachs}
\label{appendix:label}

\begin{table}[t!]
    \caption{The iteration number of CG provided with an interior optimal dual solution versus an extreme one. There are $68$, $91$, $119$ graphs in each density interval.}
    \centering
        \begin{tabular}{@{}cccc@{}}
        \toprule
         & $[0, 0.2]$ & $[0.2, 0.4]$ & $[0.4, 1.0]$ \\
    \midrule
        An interior optimal dual solution & \textbf{138.5} & \textbf{156.8} & \textbf{118.7} \\
        An extreme optimal dual solution & 325.4 & 166.4 & 123.7 \\
        \bottomrule
    \end{tabular}
    \label{tab:degen}
\end{table}

In the training data, a training example is associated with a dual variable and the label for the training example is set to the optimal solution value of the dual value. The dual problem may have multiple optimal solutions due to primal degeneracy, a variable can have different values in different optimal solutions. In this case, we set the prediction target to an average of its values in different optimal dual solutions. We first discuss the intuition behind this approach and present an empirical study to show its advantage.

Geometrically, an optimal dual problem found by simplex methods is an extreme point in the feasible dual space $P$ for the dual problem. Different optimal solutions can form a polyhedron $P^* \in P$, and any point in this polyhedron $P^*$ is another optimal dual solution. We refer to the optimal dual solution computed by an average of the extreme points to $P^*$ as an interior point. We will compare the two labeling approaches that set the target value of a dual value by its optimal solution value in an interior point and the value in an extreme point. The idea of this labeling approach is inspired by~\citep{rousseau2007interior}, which shows that using an average of multiple dual solutions to a RDP can generate better quality columns.

We find that an interior point can better boost the performance of CG compared with an extreme point, although they are both optimal dual solutions. We first identify $282$ Matilda graphs that are associated with degenerate LP instances among our selected $1326$ Matilda graphs. For a graph, we solve its LP multiple times (proportional to the number of zero-valued basis in an optimal LP solution with a factor of $10$) to collect different optimal dual solutions. Since we use the Gurobi LP solver which is based on simplex methods, these optimal dual solutions are extreme points in $P^*$. We then obtain the interior point of $P^*$ by averaging the optimal dual values in these optimal dual solutions. 

We design a method that builds a warm start using these two types of optimal dual solutions for CG. The method involves using a heuristic pricing method~\citep{shen2022enhancing} to generate $10$ columns that have the smallest reduced costs with respect to an optimal dual solution. These columns are used to initialize the RDP. The results in Table \ref{tab:degen} show that CG with the interior optimal dual solution generally uses fewer iterations to solve graphs in different density intervals compared to an extreme one, especially for sparse graphs.

\subsection{Hyperparameters and Training Procedure of ML models}
\label{appendix:mlsetup}

\begin{table*}[t!]
    \centering
    \caption{The performance of ML variants for dual solution prediction on Matilda. `Prediction time' is a summation over all test graphs. FFNN with 3 layers and GCN with 20 layers are the ML models used in our paper.}
    \label{tab:matilda_ml}
    \begin{tabular}{llcccc}
    \toprule
    \multirow{2}{*}{Method/model} & \multirow{2}{*}{\begin{tabular}{@{}l@{}} Layer/Kernel\end{tabular}} & \multicolumn{3}{c}{Mean squared error ($\times 10^{-3}$)}  & \multirow{2}{*}{\begin{tabular}{@{}c@{}} Prediction Time \\ (seconds) \end{tabular}} \\
                         &  & Training & Testing & Difference\\
    \cmidrule(lr){1-2}\cmidrule(lr){3-5}\cmidrule(lr){6-6}
                             Degree Attribute & & 4.10 & 4.53 & -0.43 & 0.0\\
    \cmidrule(lr){1-2}\cmidrule(lr){3-5}\cmidrule(lr){6-6}
    
    Logistic Regression & & 3.34 & 3.68 & -0.34 & 0.1\\
     \cmidrule(lr){1-2}\cmidrule(lr){3-5}\cmidrule(lr){6-6}

    \multirow{3}{*}{Support Vector Regression} & Linear & 3.38 & 3.75 & -0.37 & 0.1\\
                         & Polynomial & 3.50 & 3.94 & -0.44 & 3.8\\
                         & Radial Basis & 3.43 & 3.94 & -0.51 & 7.1 \\
        \cmidrule(lr){1-2}\cmidrule(lr){3-5}\cmidrule(lr){6-6}

     \multirow{3}{*}{FFNN} & 2 Layers & 3.16 & 3.51 & -0.35 & 0.2\\
                           & \textbf{3 Layers} & 2.86 & \textbf{3.34} & -0.48 & 0.2\\
                           & 5 Layers & 2.89 & 3.41 & -0.52 & 0.3\\
        \cmidrule(lr){1-2}\cmidrule(lr){3-5}\cmidrule(lr){6-6}
        \cmidrule(lr){1-2}\cmidrule(lr){3-5}\cmidrule(lr){6-6}

      \multirow{3}{*}{GCN} &  5  Layers  & 3.03 & 3.82 & -0.79 & 1.8 \\
               & 10 layers & 2.59 & 3.48 & -0.89 & 2.1 \\
               & \textbf{20 Layers}  & 1.77 & \textbf{2.99} & -1.22 & 2.5 \\
     \bottomrule
    \end{tabular}
\end{table*}

Our ML models are configured as follows. For \textit{FFNN}, we set the number of neurons for the hidden layers to $32$, and set the layer number $l$ to $3$, selected from $l \in \{2,3,5\}$. We set the parameters for GCN, following related studies~\citep{li2018combinatorial, shen2021learning}. The dimensionality of the weight matrices $\Theta^{l}$ is set to $32$, and the layer number $l$ is set to $20$, selected from the $l\in\{5,10,20\}$ layers. Besides these ML models, we have also tested other ML models including Logistic Regression and Support Vector Machine. The results can be found in Table~\ref{tab:matilda_ml}. We have chosen to present the results of our best models in the main paper for readability.


Both ML models are trained to minimize the mean squared error between predicted values and the optimal dual values for dual variables. We use the stochastic gradient descent method with Adam optimizer~\citep{kingma2014adam} to train ML models up to $1000$ epochs. In an iteration, the model parameters are updated for each training problem instance with a learning rate of~$0.0001$. To mitigate the problem of overfitting, we adopt an early stopping criterion: Stop training when an ML model cannot achieve a better validation loss for $100$ iterations. 

\subsection{Tuning the Penalty Value for Compared Methods}
\label{appendix:epstune}
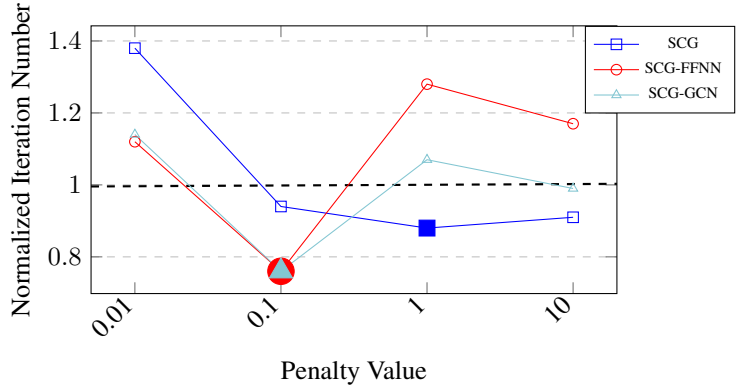
\begin{figure}[tb!]
    \centering
    
\begin{tikzpicture}
\begin{axis}[
    xlabel={Penalty Value},
    ylabel={Normalized Iteration Number},
    width=0.5\textwidth,
    height=0.3\textwidth,
    xtick=data,
    symbolic x coords={p,0.01,0.1,1,10},
    xticklabel style={rotate=45,anchor=east},
    legend style={legend columns=1,at={(1.225,1)},anchor=north east},
    ymajorgrids=true,
    grid style=dashed,
     legend style={
        font=\tiny, 
    },
]

\draw [dashed, thick] (rel axis cs:0,0.4) -- (rel axis cs:1,0.41);

\addplot[
    color=blue,
    mark=square,
]
coordinates {
    (0.01,1.38)(0.1,0.94)(1,0.88)(10,0.91)
};
\addlegendentry{SCG}

\addplot[
    color=red,
    mark=o,
]
coordinates {
    (0.01,1.12)(0.1,0.76)(1,1.28)(10,1.17)
};
\addlegendentry{SCG-FFNN}

\addplot[
    color=green,
    mark=triangle,
]
coordinates {
    (0.01,1.14)(0.1,0.76)(1,1.07)(10,0.99)
};
\addlegendentry{SCG-GCN}

\addplot[
    only marks,
    color=blue,
    mark=square*,
    mark size=3pt,
] coordinates {(1,0.88)};

\addplot[
    only marks,
    color=red,
    mark=*,
    mark size=5pt,
] coordinates {(0.1,0.76)};

\addplot[
    only marks,
    color=green,
    mark=triangle*,
    mark size=5pt,
] coordinates {(0.1,0.76)};

\node [anchor=west] at (rel axis cs:1,0.5) {CG};

\end{axis}
\end{tikzpicture}

    \caption{Performance of the compared methods with different penalty values $\epsilon$. The $x$-axis shows the iteration number normalized with respect to the iteration number of CG.}
    \label{fig:tune}
\end{figure}

Figure~\ref{fig:tune} shows the impact of setting different values for the penalty parameter for the compared methods, SCG~\citep{du1999stabilized} and SCG-FFNN/GCN~\citep{kraul2023machine}. SCG obtains its best performance with $\epsilon=1$, and SCG-FFNN/GCN obtains its best performance with $\epsilon=0.1$. Therefore, we have used these best penalty values for the compared methods. Note that the performance curve of SCG-FFNN/GCN is not `convex. This indicates that using a constant value for penalty can be insufficient and necessitates our design of the adaptive stabilization method for better utilizing a ML prediction.

\section{Extended Results for Graph Coloring Benchmarks}
\label{appendix:ret}

For reference, Table~\ref{table:gcb-class-iter} and Table~\ref{table:gcb-class-time} show the averaged iteration number and computation time for the methods with respect to graph classes. Table~\ref{table:gcb-individual-iter} and Table~\ref{table:gcb-individual-time} show the averaged iteration number and computation time for the compared methods for individual graphs.

\begin{table}[tb!]
    \centering
    \caption{Iteration number on graphs solved by all methods, averaged with respect to graph class using geometric meaning.}
    \label{table:gcb-class-iter}
    \resizebox{\textwidth}{!}{\begin{tabular}{@{}ccccccc|ccccccc@{}}
        \toprule
         \multirow{2}{*}{Instname}  & \multicolumn{6}{c}{ Exact Pricing} &  \multicolumn{6}{c}{Heuristic Pricing} \\
         
        \cmidrule(lr){2-7}\cmidrule(lr){8-13}
             & CG & SCG & SCG-FFNN & SCG-GCN & ASCG-FFNN & ASCG-GCN & CG & SCG & SCG-FFNN & SCG-GCN & ASCG-FFNN & ASCG-GCN \\
                                     
        \cmidrule(lr){1-1}\cmidrule(lr){2-7}\cmidrule(lr){8-13}
        DSJC(6) & 560.4 & 507.1 & 379.9 & 402.4 & 349.5 & 368.7 & 692.7 & 725.9 & 518.4 & 526.7 & 505.9 & 522.8 &  \\
        FullIns(6) & 63.2 & 105.2 & 65.1 & 86.7 & 31.9 & 58.6 & 104.1 & 114.9 & 90.5 & 101.1 & 78.9 & 98.1 &  \\
        Insertions(4) & 358.3 & 176.4 & 150.5 & 281.5 & 133.5 & 269.1 & 515.7 & 506.1 & 468.1 & 476.4 & 404.6 & 492.5 &  \\
        flat(1) & 1240.9 & 816.4 & 832.8 & 975.5 & 770.1 & 954.2 & 1812.7 & 1407.2 & 1290.9 & 1493.2 & 1261.4 & 1487.8 &  \\
        le(0) & nan & nan & nan & nan & nan & nan & nan & nan & nan & nan & nan & nan &  \\
        myciel(4) & 71.8 & 60.9 & 48.2 & 55.4 & 38.2 & 49.9 & 100.7 & 112.0 & 93.4 & 97.4 & 77.7 & 90.7 &  \\
        queen(9) & 91.5 & 158.0 & 73.0 & 92.6 & 61.0 & 80.5 & 120.8 & 166.5 & 92.0 & 109.2 & 81.1 & 98.1 &  \\
        \bottomrule
    \end{tabular}}
\end{table}

\begin{table}[tb!]
    \centering
    \caption{Computation time for graphs solved by at least one method, averaged with respect to graph class using geometric meaning.}
    \label{table:gcb-class-time}
    \resizebox{\textwidth}{!}{\begin{tabular}{@{}ccccccc|ccccccc@{}}
        \toprule
         \multirow{2}{*}{Instname}  & \multicolumn{6}{c}{ Exact Pricing} &  \multicolumn{6}{c}{Heuristic Pricing} \\
         
        \cmidrule(lr){2-7}\cmidrule(lr){8-13}
             & CG & SCG & SCG-FFNN & SCG-GCN & ASCG-FFNN & ASCG-GCN & CG & SCG & SCG-FFNN & SCG-GCN & ASCG-FFNN & ASCG-GCN \\
                                     
        \cmidrule(lr){1-1}\cmidrule(lr){2-7}\cmidrule(lr){8-13}
        DSJC(6) & 28.1 & 33.5 & 30.4 & 30.2 & 25.9 & 26.5 & 28.5 & 36.6 & 27.8 & 28.9 & 25.4 & 25.7 &  \\
        FullIns(9) & 6.8 & 20.7 & 14.9 & 18.1 & 10.9 & 10.2 & 13.3 & 13.8 & 12.7 & 13.3 & 12.4 & 12.5 &  \\
        Insertions(6) & 33.0 & 12.0 & 18.1 & 30.3 & 16.6 & 27.3 & 39.1 & 39.3 & 38.2 & 38.2 & 35.9 & 39.3 &  \\
        flat(1) & 514.3 & 432.3 & 454.8 & 495.4 & 433.7 & 495.7 & 473.4 & 428.4 & 415.6 & 465.7 & 407.0 & 465.0 &  \\
        le(2) & 48.2 & 1000.0 & 87.8 & 99.3 & 91.5 & 84.7 & 76.2 & 77.6 & 145.6 & 76.1 & 85.3 & 81.6 &  \\
        myciel(5) & 7.0 & 3.9 & 5.8 & 6.2 & 5.5 & 5.9 & 7.5 & 7.8 & 7.5 & 7.5 & 7.1 & 7.3 &  \\
        queen(13) & 19.1 & 34.7 & 16.4 & 20.1 & 14.5 & 19.2 & 10.4 & 9.4 & 7.5 & 9.8 & 7.6 & 9.8 &  \\
        \bottomrule
    \end{tabular}}
\end{table}

\begin{table}[tb!]
    \centering
    \caption{Iteration number on individual graphs solved by all methods.}
    \label{table:gcb-individual-iter}
    \resizebox{\textwidth}{!}{\begin{tabular}{@{}rrrrrrr|rrrrrrr@{}}
        \toprule
         \multirow{2}{*}{Instname}  & \multicolumn{6}{c}{ Exact Pricing} &  \multicolumn{6}{c}{Heuristic Pricing} \\
         
        \cmidrule(lr){2-7}\cmidrule(lr){8-13}
             & CG & SCG & SCG-FFNN & SCG-GCN & ASCG-FFNN & ASCG-GCN & CG & SCG & SCG-FFNN & SCG-GCN & ASCG-FFNN & ASCG-GCN \\
                                     
        \cmidrule(lr){1-1}\cmidrule(lr){2-7}\cmidrule(lr){8-13}
        
        DSJC125.1 & 1132.8 & 544.3 & 577.5 & 788.3 & 518.6 & 676.8 & 1755.5 & 1729.0 & 1671.8 & 1727.3 & 1715.0 & 1678.7\\
        DSJC125.5 & 410.0 & 331.9 & 233.3 & 250.7 & 215.4 & 218.9 & 515.2 & 414.1 & 312.0 & 322.7 & 294.7 & 310.0\\
        DSJC125.9 & 167.8 & 241.3 & 143.5 & 136.0 & 136.9 & 130.4 & 187.2 & 271.0 & 143.0 & 137.3 & 140.6 & 137.2\\
        DSJC250.5 & 976.7 & 673.3 & 641.8 & 681.8 & 573.9 & 629.3 & 1334.1 & 1050.0 & 924.7 & 983.5 & 875.8 & 949.5\\
        DSJC250.9 & 431.8 & 518.5 & 315.5 & 308.7 & 292.7 & 287.8 & 466.7 & 565.5 & 334.1 & 336.8 & 323.3 & 330.5\\
        DSJC500.9 & 939.5 & 1115.8 & 765.0 & 748.3 & 706.6 & 714.9 & 1044.9 & 1267.6 & 838.5 & 838.2 & 829.7 & 906.9\\
        \cmidrule(lr){1-1}\cmidrule(lr){2-7}\cmidrule(lr){8-13}
        1-FullIns\_3 & 41.0 & 58.0 & 29.0 & 73.0 & 13.0 & 54.0 & 61.2 & 67.5 & 34.8 & 52.1 & 19.7 & 52.9\\
        1-FullIns\_4 & 422.0 & 297.0 & 684.0 & 930.0 & 148.0 & 260.0 & 967.8 & 1048.6 & 887.4 & 972.8 & 879.1 & 875.0\\
        2-FullIns\_3 & 47.0 & 105.0 & 62.0 & 58.0 & 22.0 & 57.0 & 59.6 & 64.3 & 55.6 & 61.5 & 48.2 & 61.0\\
        3-FullIns\_3 & 30.0 & 168.0 & 32.0 & 44.0 & 24.0 & 32.0 & 79.4 & 87.3 & 72.4 & 78.3 & 64.6 & 73.4\\
        4-FullIns\_3 & 53.0 & 64.0 & 39.0 & 42.0 & 29.0 & 37.0 & 56.9 & 65.6 & 55.4 & 56.4 & 50.1 & 53.9\\
        5-FullIns\_3 & 48.0 & 69.0 & 48.0 & 57.0 & 34.0 & 42.0 & 78.5 & 86.8 & 77.8 & 76.2 & 85.7 & 78.0\\
        \cmidrule(lr){1-1}\cmidrule(lr){2-7}\cmidrule(lr){8-13}
        1-Insertions\_4 & 490.0 & 182.0 & 156.0 & 454.0 & 129.0 & 488.0 & 800.5 & 812.1 & 797.9 & 793.8 & 745.3 & 806.4\\
        2-Insertions\_3 & 142.0 & 101.0 & 76.0 & 101.0 & 64.0 & 89.0 & 162.9 & 159.9 & 120.9 & 136.6 & 85.8 & 154.2\\
        3-Insertions\_3 & 348.0 & 177.0 & 137.0 & 223.0 & 127.0 & 202.0 & 469.4 & 444.4 & 450.1 & 426.4 & 386.6 & 433.3\\
        4-Insertions\_3 & 679.0 & 297.0 & 315.0 & 612.0 & 302.0 & 595.0 & 1152.3 & 1133.6 & 1101.3 & 1110.1 & 1076.8 & 1088.6\\
        \cmidrule(lr){1-1}\cmidrule(lr){2-7}\cmidrule(lr){8-13}
        flat300\_28\_0 & 1240.9 & 816.4 & 832.8 & 975.5 & 770.1 & 954.2 & 1812.7 & 1407.2 & 1290.9 & 1493.2 & 1261.4 & 1487.8\\
        \cmidrule(lr){1-1}\cmidrule(lr){2-7}\cmidrule(lr){8-13}
        myciel3 & 12.0 & 17.0 & 14.0 & 14.0 & 10.0 & 11.0 & 12.8 & 19.4 & 13.0 & 14.6 & 9.8 & 10.9\\
        myciel4 & 33.0 & 39.0 & 30.0 & 31.0 & 22.0 & 28.0 & 48.7 & 46.8 & 36.9 & 39.9 & 26.9 & 41.2\\
        myciel5 & 128.0 & 91.0 & 61.0 & 68.0 & 47.0 & 65.0 & 206.1 & 206.4 & 190.4 & 192.8 & 166.2 & 188.8\\
        myciel6 & 491.0 & 220.0 & 202.0 & 304.0 & 194.0 & 291.0 & 752.2 & 804.5 & 780.3 & 757.6 & 762.1 & 739.5\\
        \cmidrule(lr){1-1}\cmidrule(lr){2-7}\cmidrule(lr){8-13}
        queen5\_5 & 8.1 & 19.7 & 11.0 & 9.0 & 7.5 & 5.5 & 12.7 & 23.3 & 9.0 & 9.0 & 6.0 & 5.8\\
        queen6\_6 & 62.5 & 95.1 & 42.4 & 44.1 & 39.6 & 41.5 & 76.2 & 100.7 & 42.8 & 45.2 & 39.1 & 42.1\\
        queen7\_7 & 14.0 & 58.8 & 26.7 & 33.6 & 21.9 & 22.1 & 15.1 & 61.0 & 16.1 & 15.6 & 10.1 & 11.6\\
        queen8\_12 & 17.3 & 272.6 & 87.6 & 90.3 & 84.4 & 84.7 & 19.6 & 220.4 & 86.3 & 97.3 & 86.0 & 93.3\\
        queen8\_8 & 148.8 & 182.8 & 87.3 & 95.4 & 68.8 & 85.0 & 197.2 & 190.8 & 103.8 & 112.4 & 79.3 & 101.1\\
        queen9\_9 & 281.7 & 192.0 & 104.4 & 164.4 & 89.6 & 147.4 & 389.0 & 277.7 & 190.7 & 258.7 & 215.3 & 260.5\\
        queen10\_10 & 380.0 & 292.1 & 193.8 & 234.9 & 169.0 & 221.7 & 514.0 & 346.7 & 316.3 & 386.9 & 296.4 & 374.2\\
        queen11\_11 & 407.2 & 366.7 & 154.6 & 298.4 & 123.9 & 295.5 & 591.3 & 379.4 & 323.3 & 480.3 & 316.7 & 465.8\\
        queen12\_12 & 473.4 & 521.9 & 184.4 & 342.1 & 147.8 & 343.3 & 712.3 & 431.1 & 385.4 & 574.9 & 378.7 & 567.5\\
        \bottomrule
    \end{tabular}}
\end{table}

\begin{table}[tb!]
    \centering
    \caption{Computation time for individual graphs solved by at least one method.}
    \label{table:gcb-individual-time}
    \resizebox{\textwidth}{!}{\begin{tabular}{@{}rrrrrrr|rrrrrrr@{}}
        \toprule
         \multirow{2}{*}{Instname}  & \multicolumn{6}{c}{ Exact Pricing} &  \multicolumn{6}{c}{Heuristic Pricing} \\
         
        \cmidrule(lr){2-7}\cmidrule(lr){8-13}
             & CG & SCG & SCG-FFNN & SCG-GCN & ASCG-FFNN & ASCG-GCN & CG & SCG & SCG-FFNN & SCG-GCN & ASCG-FFNN & ASCG-GCN \\
                                     
        \cmidrule(lr){1-1}\cmidrule(lr){2-7}\cmidrule(lr){8-13}
        DSJC125.1 & 176.2 & 104.8 & 209.2 & 163.3 & 139.6 & 133.8 & 210.8 & 212.6 & 204.3 & 213.0 & 204.1 & 203.2 \\
        DSJC125.5 & 8.6 & 8.7 & 5.8 & 6.6 & 5.3 & 5.6 & 8.7 & 8.4 & 5.7 & 6.4 & 5.4 & 5.9 \\
        DSJC125.9 & 1.5 & 2.5 & 1.5 & 1.4 & 1.4 & 1.3 & 1.2 & 2.1 & 1.1 & 1.0 & 1.0 & 1.0 \\
        DSJC250.5 & 153.1 & 136.1 & 134.4 & 140.2 & 119.7 & 128.8 & 145.4 & 138.1 & 123.6 & 134.4 & 108.0 & 119.8 \\
        DSJC250.9 & 9.2 & 15.8 & 11.6 & 11.2 & 10.0 & 9.9 & 8.5 & 15.2 & 10.6 & 10.7 & 9.3 & 9.2 \\
        DSJC500.9 & 90.3 & 201.2 & 158.6 & 178.9 & 130.9 & 148.0 & 106.9 & 200.4 & 140.3 & 140.6 & 112.2 & 101.2 \\
        \cmidrule(lr){1-1}\cmidrule(lr){2-7}\cmidrule(lr){8-13}
        1-FullIns\_3 & 0.3 & 0.4 & 0.2 & 0.5 & 0.1 & 0.4 & 0.2 & 0.3 & 0.1 & 0.2 & 0.1 & 0.2 \\
        1-FullIns\_4 & 9.1 & 24.5 & 28.5 & 63.1 & 3.6 & 8.1 & 21.2 & 26.0 & 18.8 & 21.6 & 16.7 & 17.5 \\
        2-FullIns\_3 & 0.3 & 1.1 & 0.6 & 0.5 & 0.2 & 0.5 & 0.3 & 0.3 & 0.3 & 0.3 & 0.2 & 0.3 \\
        3-FullIns\_3 & 0.2 & 3.1 & 0.3 & 0.4 & 0.3 & 0.3 & 0.4 & 0.5 & 0.4 & 0.4 & 0.4 & 0.4 \\
        4-FullIns\_3 & 0.5 & 0.7 & 0.5 & 0.4 & 0.6 & 0.4 & 0.4 & 0.4 & 0.4 & 0.4 & 0.4 & 0.4 \\
        5-FullIns\_3 & 0.6 & 1.0 & 0.8 & 0.8 & 0.8 & 0.6 & 0.7 & 0.8 & 0.7 & 0.7 & 0.8 & 0.7 \\
        2-FullIns\_4 & 72.8 & Cutoff & 327.3 & 686.1 & 198.4 & 308.3 & 235.7 & 226.2 & 193.8 & 239.5 & 210.0 & 188.7 \\
        3-FullIns\_4 & 116.0 & Cutoff & Cutoff & 932.9 & Cutoff & 157.3 & 877.2 & 915.6 & 899.7 & 885.5 & 879.3 & 843.7 \\
        4-FullIns\_4 & 243.5 & Cutoff & Cutoff & Cutoff & Cutoff & Cutoff & Cutoff & Cutoff & Cutoff & Cutoff & Cutoff & Cutoff \\
        \cmidrule(lr){1-1}\cmidrule(lr){2-7}\cmidrule(lr){8-13}
        1-Insertions\_4 & 8.4 & 2.6 & 2.1 & 9.2 & 1.7 & 9.8 & 12.2 & 12.8 & 12.6 & 12.5 & 10.9 & 12.5 \\
        2-Insertions\_3 & 1.3 & 1.0 & 0.7 & 1.0 & 0.5 & 0.8 & 0.7 & 0.8 & 0.6 & 0.7 & 0.6 & 1.0 \\
        3-Insertions\_3 & 4.6 & 2.1 & 1.6 & 3.1 & 1.5 & 2.7 & 4.6 & 4.4 & 4.4 & 4.2 & 3.6 & 4.2 \\
        4-Insertions\_3 & 19.5 & 6.7 & 9.4 & 20.4 & 8.8 & 18.9 & 31.8 & 31.7 & 30.3 & 30.5 & 28.5 & 28.9 \\
        1-Insertions\_5 & Cutoff & 236.2 & Cutoff & Cutoff & Cutoff & Cutoff & Cutoff & Cutoff & Cutoff & Cutoff & Cutoff & Cutoff \\
        2-Insertions\_4 & 623.1 & 121.6 & 336.4 & 524.3 & 291.6 & 361.2 & 986.5 & 998.8 & 999.4 & 986.6 & Cutoff & Cutoff \\
        \cmidrule(lr){1-1}\cmidrule(lr){2-7}\cmidrule(lr){8-13}
        flat300\_28\_0 & 514.3 & 432.3 & 454.8 & 495.4 & 433.7 & 495.7 & 473.4 & 428.4 & 415.6 & 465.7 & 407.0 & 465.0 \\
        \cmidrule(lr){1-1}\cmidrule(lr){2-7}\cmidrule(lr){8-13}
        myciel3 & 0.1 & 0.1 & 0.1 & 0.1 & 0.1 & 0.1 & 0.1 & 0.1 & 0.1 & 0.1 & 0.0 & 0.0 \\
        myciel4 & 0.2 & 0.3 & 0.2 & 0.2 & 0.1 & 0.2 & 0.2 & 0.2 & 0.2 & 0.2 & 0.1 & 0.2 \\
        myciel5 & 1.1 & 0.9 & 0.6 & 0.7 & 0.4 & 0.6 & 1.3 & 1.4 & 1.2 & 1.3 & 1.0 & 1.2 \\
        myciel6 & 10.7 & 5.5 & 5.6 & 7.4 & 5.7 & 7.1 & 13.9 & 15.8 & 14.8 & 14.5 & 14.1 & 13.7 \\
        myciel7 & Cutoff & 152.9 & Cutoff & Cutoff & Cutoff & Cutoff & Cutoff & Cutoff & Cutoff & Cutoff & Cutoff & Cutoff \\
        \cmidrule(lr){1-1}\cmidrule(lr){2-7}\cmidrule(lr){8-13}
        queen5\_5 & 0.1 & 0.1 & 0.1 & 0.1 & 0.0 & 0.0 & 0.0 & 0.1 & 0.0 & 0.0 & 0.0 & 0.0 \\
        queen6\_6 & 0.4 & 0.7 & 0.3 & 0.3 & 0.3 & 0.3 & 0.3 & 0.4 & 0.2 & 0.2 & 0.2 & 0.2 \\
        queen7\_7 & 0.1 & 0.5 & 0.2 & 0.3 & 0.2 & 0.2 & 0.1 & 0.3 & 0.1 & 0.1 & 0.0 & 0.0 \\
        queen8\_12 & 0.2 & 12.2 & 2.5 & 2.2 & 2.3 & 2.0 & 0.1 & 2.9 & 0.7 & 0.9 & 0.8 & 0.8 \\
        queen8\_8 & 1.5 & 2.1 & 1.0 & 1.1 & 0.7 & 0.9 & 1.1 & 1.2 & 0.7 & 0.8 & 0.5 & 0.7 \\
        queen9\_9 & 4.2 & 2.9 & 1.6 & 2.7 & 1.3 & 2.3 & 3.8 & 2.6 & 1.9 & 2.9 & 2.0 & 2.7 \\
        queen10\_10 & 9.9 & 10.0 & 5.9 & 7.1 & 4.9 & 6.5 & 7.7 & 4.4 & 4.8 & 6.5 & 4.2 & 5.7 \\
        queen11\_11 & 19.8 & 56.2 & 9.4 & 20.8 & 7.4 & 21.4 & 12.4 & 7.4 & 6.7 & 11.8 & 6.0 & 10.5 \\
        queen12\_12 & 63.5 & 517.8 & 37.5 & 65.6 & 27.6 & 71.4 & 22.4 & 15.1 & 12.2 & 21.0 & 11.4 & 19.3 \\
        queen13\_13 & 242.2 & Cutoff & 158.8 & 295.3 & 135.8 & 261.8 & 37.0 & 23.3 & 18.7 & 33.3 & 17.9 & 32.4 \\
        queen14\_14 & Cutoff & Cutoff & 940.7 & Cutoff & 752.2 & Cutoff & 100.2 & 69.5 & 53.1 & 89.6 & 58.7 & 89.4 \\
        queen15\_15 & Cutoff & Cutoff & Cutoff & Cutoff & Cutoff & Cutoff & 444.0 & 303.6 & 206.9 & 302.4 & 324.1 & 420.1 \\
        queen16\_16 & Cutoff & Cutoff & Cutoff & Cutoff & Cutoff & Cutoff & 749.0 & 680.3 & 887.7 & 781.1 & 911.0 & 915.9 \\
        \cmidrule(lr){1-1}\cmidrule(lr){2-7}\cmidrule(lr){8-13}
        le450\_25a & 49.1 & Cutoff & 94.5 & 114.0 & 103.7 & 87.9 & 88.3 & 95.4 & 141.0 & 87.1 & 92.1 & 101.4 \\
        le450\_25b & 47.3 & Cutoff & 81.5 & 86.4 & 80.6 & 81.6 & 65.8 & 63.1 & 150.4 & 66.4 & 78.9 & 65.6 \\
        \bottomrule
    \end{tabular}}
\end{table}
\end{document}